\documentclass[11pt]{amsart}
\usepackage{amsmath}
\usepackage{amsfonts}
\usepackage{amssymb}
\usepackage{amsthm}
\usepackage{url}
\usepackage[all]{xy}
\usepackage{geometry}
\usepackage{dsfont}
\usepackage{graphicx}
\usepackage{caption}
\usepackage{subcaption}
\usepackage{comment} 
\usepackage{stmaryrd}
\usepackage{hyperref}
\usepackage{todonotes}
\usepackage{color}
\usepackage{enumerate,tikz-cd}
\usepackage{fixltx2e}
\usepackage{MnSymbol}
\usepackage{comment}
\usepackage{mathtools}
\usepackage[shortlabels]{enumitem}
\allowdisplaybreaks

\numberwithin{equation}{section}

\newtheorem{proposition}{Proposition}[section]
\newtheorem{lemma}[proposition]{Lemma}
\newtheorem{theorem}[proposition]{Theorem}
\newtheorem{corollary}[proposition]{Corollary}

\theoremstyle{definition}
\newtheorem{remark}[proposition]{Remark}
\newtheorem{definition}[proposition]{Definition}

\DeclareMathOperator{\ord}{ord}
\DeclareMathOperator{\QM}{QM}
\newcommand{\Xlin}{X^{\mathrm{lin}}}

\DeclareMathOperator{\rees}{Rees}
\DeclareMathOperator{\Bigc}{Big}
\DeclareMathOperator{\de}{d_1}
\DeclareMathOperator{\dd}{d}
\DeclareMathOperator{\nn}{N}
\DeclareMathOperator{\qd}{\underline{\mathrm{d}}}
\newcommand{\Xval}{X^{\mathrm{val}}}
\newcommand{\Xdiv}{X^{\mathrm{div}}}
\newcommand{\Xqm}{X^{\mathrm{qm}}}
\newcommand{\Xna}{X^{\mathrm{na}}}

\newcommand{\N}{\mathbb{N}}
\newcommand{\R}{\mathbb{R}}

\newcommand{\Z}{\mathbb{Z}}

\newcommand{\Q}{\mathbb{Q}}

\renewcommand{\epsilon}{\varepsilon}

\renewcommand{\phi}{\varphi}

\newcommand\MA{\mathrm{MA}}
\newcommand\FS{\mathrm{FS}}

\newcommand\mi{^{-1}}
\newcommand\vol{\mathrm{vol}}
\newcommand\na{{\mathrm{NA}}}

\newcommand\bbp{\mathbb{P}}
\newcommand\bbq{\mathbb{Q}}
\newcommand\bbr{\mathbb{R}}

\newcommand\cD{\mathcal{D}}

\newcommand{\cF}{\mathcal{F}}

\newcommand\cL{\mathcal{L}}
\newcommand\cM{\mathcal{M}}
\newcommand\cN{\mathcal{N}}
\newcommand\cO{\mathcal{O}}

\newcommand\cT{\mathcal{T}}

\newcommand\cX{\mathcal{X}}

\newcommand{\norm}[1][\cdot]{\left\|#1\right\|}

\newcommand\fa{\mathfrak{a}}

\newcommand\triv{{\mathrm{triv}}}

\renewcommand\hom{{\mathrm{hom}}}

\renewcommand{\=}{\coloneqq}
\newcommand\qm{{\mathrm{qm}}}
\newcommand{\lin}{\mathrm{lin}}

\renewcommand\div{\mathrm{div}}

\renewcommand\hom{\mathrm{hom}}

\title[Quasimonomial Monge--Ampère equations]{Quasimonomial Monge--Ampère equations for big classes}
\author[R\'emi Reboulet and Yueqiao Wu]{R\'emi Reboulet and Yueqiao Wu}

\address{R\'emi Reboulet, CNRS, Institut Camille Jordan, 21 Av. Claude Bernard, 69100 Villeurbanne, France}\email{reboulet@math.univ-lyon1.fr}
\address{Yueqiao Wu, Department of Mathematics, Johns Hopkins University, 3400 N. Charles Street
Baltimore MD 21218, United States}\email{ywu347@jh.edu}

\begin{document}

\begin{abstract}
Let $X$ be a normal projective variety over an algebraically closed field of characteristic zero. We describe an alternative approach to a differentiability result of volumes of norms defined by finitely many quasimonomial valuations. As a consequence, we prove the non-Archimedean Calabi--Yau theorem for quasimonomial measures on $X$ with respect to big classes, in the trivially valued case. We also obtain continuity of solutions as numerical classes vary in the big cone. This extends previous results of Boucksom--Jonsson and Dervan--Reboulet, and gives the first uniqueness results in the big case. Our approach relies on convex geometry and Okounkov bodies, without involving non-Archimedean pluripotential theory or test configurations.
\end{abstract}
\maketitle

\setcounter{tocdepth}{1}
\tableofcontents

\section*{Introduction}

The Yau--Tian--Donaldson (or \textit{YTD}) conjecture asks for an algebraic condition on a pair $(X,L)$ of a complex projective variety with an ample line bundle, that is necessary and sufficient for existence of a constant scalar curvature Kähler metric on $L$.

This algebraic condition, now known as \textit{K-stability}, was initially formulated by Donaldson \cite{donaldson:scalar} in terms of the positivity of a certain invariant along all \textit{test configurations} for $(X,L)$, which are certain $\mathbb{G}_m$-equivariant degenerations of $(X,L)$. This point of view, informed by geometric invariant theory, has been simplified to a more valuative approach, first pioneered by the fundamental work of Li--Xu \cite{lixu} in the Fano setting (see also \cite{dar:zhangbigke,derreb:1} in the case $-K_X$ big). In the general ample case, work of Dervan--Legendre \cite{dervanlegendre} and then Boucksom--Jonsson \cite{bj:kstab1,bj:kstab2} have led to the formulation of the notion of \textit{divisorial stability}, which tests so-called \textit{divisorial measures}, convex weighted finite sets of divisorial valuations on $X$.

The notion of divisorial stability has been particularly fertile, giving the adequate condition to a  resolution of the YTD conjecture in the ample case \cite{bj:ytd} (see also \cite{darzhang:ytd,trusiani:fujita}). There are two essential aspects to working with sets of divisorial measures: this class of test objects is independent of the choice of polarisation (so that one can study how objects involved in divisorial stability vary as the line bundle varies, which was crucial in the proof of openness in \cite{bj:kstab2}, and also allows to extend the theory to big line bundles, as was done by the first author in \cite{derreb:3}); and independent of the birational equivalence class of $X$, a fact which was used in the aforementioned \cite{derreb:3} to prove birational invariance of divisorial stability.

The heart of the approach of \cite{bj:kstab2} to define divisorial stability is the \textit{non-Archimedean Calabi--Yau theorem}. Namely, to any ample test configuration $(\cX,\cL)$ for $(X,L)$, $L$ ample, one associates a divisorial measure as before, its non-Archimedean Monge--Ampère measure, given by intersection numbers on the central fibre. This operation can be abstractly extended to more general filtrations by density, and the Calabi--Yau theorem (in this case) states that any divisorial measure is the Monge--Ampère measure of a divisorial filtration of the section ring of $L$. In \cite{derreb:3} this specific result is extended to the big line bundle case as well, analogously to the complex Archimedean case \cite{begz,bbgz}.

Let $X$ be a normal projective variety over a trivially valued algebraically closed field of characteristic zero. Our goal of this paper is to describe an alternative approach to non-Archimedean Monge--Amp\`ere equations for measures supported on a finite set of quasimononmial valuations, with respect to big classes. 

\subsection*{Divisorial Monge--Amp\`ere equations for polarized projective varieties}
Let us briefly recall the theory of non-Archimedean Monge--Amp\`ere equations for a polarized projective variety $(X, L)$, for $L$ ample, developed in~\cite{bj:trivval}, which in turn builds on~\cite{bfj:semi, bfj:nama, bhj:duistermaat}.

Through a Rees correspondence, ample normal test configurations for $(X, L)$ are in one-to-one correspondence with (decreasing, left-continuous, multiplicative) finitely generated filtrations on the section ring $R(X,L)$, which further correspond to finitely generated $k^\times$-invariant (homogeneous, non-Archimedean) norms $\chi\colon R(X, L)\to\Q\cup \{\infty\}$. The space $\cT_\Q(L)$ of norms $\chi_{\cX, \cL}$ associated to normal ample test configurations sits as a natural subspace in the space $\cN^\div(L)$ of divisorial norms on $R(X,L)$, consisting of norms of the form $\chi=\min_{v\in\Sigma}\{v+t_v\}$ for $\Sigma\subset \Xdiv$ a finite set of divisorial valuations, and $t_v\in \R$. In particular, the norm $\chi_{\cX, \cL}$ associated to a test configuration is given by
\[
\chi_{\cX,\cL}=\min_{E\subseteq\cX_0}\{v_E+t_E\},
\]
where $E$ runs through all irreducible components of the central fiber $\cX_0\subset\cX$ and $v_E$ is a divisorial valuation associated to $E$. The \emph{Monge--Amp\`ere operator} with respect to $L$,
\[
\MA_L\colon \cT_\Q(L)\to \cM^\div, \,\chi_{\cX, \cL}\mapsto \MA_L(\chi_{\cX,\cL})\=\vol(L)^{-1}\sum_{E\subset\cX_0}\ord_E(\cX_0)(\cL|_E)^n\delta_{v_E},
\]
assigns to any norm in $\cT_\Q$ a \emph{divisorial measure}, i.e.\ 
a probability measure on $\Xna$ supported in a finite set $\Sigma=\{v_E\mid E\subset \cX_0\}$ of divisorial valuations. 

Partly motivated by the study of metric structures on the space of K\"ahler potentials in~\cite{dar:mabuchi}, the space $\cT_\Q(L)$ and more generally $\cN_\R^\div(L)$ comes with a natural Darvas metric defined through volumes of filtrations, as studied in~\cite{bj:kstab1,reb:1}. More precisely, to each test configuration $(\cX, \cL)$, or equivalently norm $\chi_{\cX, \cL}$, is associated a (family of) graded linear series $(V^\lambda_\bullet)_{\lambda\in\R}$ of $R(X,L)$, defined by
\[
V_m^\lambda(\chi)\= \{s\in H^0(X, mL)\mid \chi_{\cX, \cL}(s)\ge m\lambda\}.
\]
The \emph{volume} of $\chi$ is the average of the volumes of the graded linear series $V^\lambda_\bullet$:
\[
\vol(\chi)\=c+\frac{1}{\vol(L)}\int_c^\infty \vol(V^\lambda_\bullet)\, d\lambda,
\]
for $c\ll 0$. The metric $\de$ on $\cN^\div(L)$ is defined by setting, for each $\chi, \chi'\in \cN^\div(L)$,
\[
\de(\chi, \chi')=\vol(\chi)+\vol(\chi')-2\vol(\chi\wedge\chi'),
\]
where $\chi\wedge \chi'$ denotes the pointwise minimum. Equipped with the $\de$ metric, the spacd $\cT_\Q(L)$ is a dense subset of $\cN^\div(L)$ and it is shown in~\cite{bj:kstab1} that the Monge--Amp\`ere operator extends continuously to an isomorphism
\[
\MA_L\colon \cN^\div(L)/\R\to \cM^\div,
\]
where $\cN^\div(L)/\R$ is the quotient by translations equipped with the quotient metric of $\de$.


\subsection*{Differentiability of volumes of quasimonomial data}
The starting point of our approach is to consider the space $\cD(\Xqm)$ of \emph{quasimonomial data}, consisting of pairs $(\Sigma\=\{v_i\}_{i=1}^d, t)$ where $\Sigma$ is a finite set of quasimonomial valuations and $t\in \R^d$. This space has the advantage of being independent of the polarization, and we will study invariants on $\cD(\Xqm)$ as we vary numerical classes in the big cone. More precisely, for each big line bundle $L$,
the space $\cD(\Xqm)$ comes with a natural map to the space $\cN^\qm(L)$ of \emph{quasimonomial norms} on $L$, generalizing the previously discussed divisorial norms, given by
\[
\nn(L)\colon \cD(\Xqm)\to \cN^\qm(L),\quad (\Sigma, t)\mapsto \chi_{\Sigma, t}\=\min_{v\in\Sigma}\{v+t_v\},
\]
on $R(X, L)$; and the volume function
\[
\vol(L; \Sigma, \cdot)\colon \R^d\to \R,\quad t\mapsto \vol(L;\Sigma,t)\=\vol_L(\chi_{\Sigma, t})
\]
recovers the expected vanishing order of a valuation (if $d=1$), the non-Archimedean Monge--Ampère energy of a test configuration, and more generally the volume type invariants appearing in divisorial stability as in~\cite{wn:pietro,derreb:3}. Our first result states that this function only depends on the numerical class of $L$, extends continuously to the big cone and defines a differentiable function on $\R^d$.

\bigskip\noindent\textbf{Theorem A.} 
For a fixed finite subset $\Sigma\subset \Xqm$ of size $d$, we have 
\begin{itemize}
    \item [(i)] for any $t\in \R^d$, the function $\vol(\cdot; \Sigma, t)\colon \Bigc(X)\to \R$ is continuous;
    \item[(ii)] for $\theta\in\Bigc(X)$, $\vol(\theta; \Sigma, t)\in C^1(\R^d)$.
\end{itemize}


\bigskip 
While the volume of a norm on a big line bundle $L$ is not defined by a single graded linear series, it has an alternative characterization using the concave transform on an Okounkov body of $L$. The proof of Theorem~A is thus reminiscent of the classical differentiability result of the volume function on the big cone. The volume function readily extends to $\Bigc(X)_\Q$ by homogeneity, and thus extends to uniquely to a continuous function on the big cone once we have proved local Lipschitz continuity. The proof of differentiability is in turn inspired by the Okounkov body construction in~\cite{lazmus}, where they used it to prove that the classical volume function is $C^1$ on the big cone. Since we consider quasimonomial valuations in general, the Okounkov body construction we use comes from the work of Peng in~\cite{peng:stab}, who defined Okounkov bodies with respect to quasimonomial valuations. 

By nature of our proof, we obtain explicit expressions of $\nabla\vol(\theta; \Sigma, t)\=(\MA_\theta(\Sigma, t)_i)_{i=1}^d$, given by
$$\MA_\theta(\Sigma,t)_i:=\frac{n!}{\norm[\alpha_i]\vol(L)}\int_\bbr \vol(\Delta_i^\lambda(\theta;t)\cap \{x\in \R^d\mid \langle\alpha_i, x\rangle +t_i=\lambda\})\,d\lambda,$$
where $\Delta_i^\lambda(\theta;t)$ is a filtered Okounkov body with respect to the valuation $v_i$ and $\alpha_i$ is the weight associated to $v_i$.
Even in the the case when $(\Sigma,t)$ defines a divisorial norm $\chi\in \cT_\Q(L)$ for $\theta=c_1(L)$ ample, to the best of the knowledge of the authors, this specializes to give a (new) explicit expression of the Monge--Ampère measure $\MA(\chi)=\sum_i\MA(\Sigma, t)_i\delta_{v_i}$ without relying on test configurations: if $v_i=c_i^{-1}\ord_{ E_i}$, then
$$\MA_\theta(\Sigma,t)_i=\frac{n}{\vol(L)}\left(c_i\mi\int_{t_i}^\infty\vol_{Y|E_i}(\pi^*L-\sum_{j=1}^dc_j^{-1}(\lambda-t_j)_+E_j)\,d\lambda\right).$$

\subsection*{Quasimonomial Monge--Amp\`ere equations}
Theorem~A allows us to define the \textit{Monge--Ampère operator} on $\cD(\Xqm)$ as the map
\[
\MA_\theta\colon \cD(\Xqm)\to \cM^\qm, \quad (\Sigma=\{v_i\}_{i=1}^d, t)\mapsto \MA_\theta(\Sigma,t)\=\sum_{i=1}^d \MA_\theta(\Sigma,t)_i \delta_{v_i},
\]
where $\cM^\qm$ is the space of \emph{quasimonimial measures}, that is, probability measures on $\Xna$ supported on a finite set of quasimonomial valuations.
Recall that we have a natural map $\nn(L)\colon \cD(\Xqm)\to \cN^\qm(L)$ for each big line bundle $L$, which by construction is surjective but far from being injective. Thus it is not a priori clear that our notion of the Monge--Amp\`ere operator recoves the one defined by Boucksom--Jonsson in the case when $L$ is ample.

To this end, let us note that Theorem~A in particular enables us to extend the metric structures on $\cN^\qm(L)$ to a family of pseudometrics $\{\dd_{1,\theta}\}_{\theta\in\Bigc(X)}$ on $\cD(\Xqm)$, by setting for $(\Sigma, t), (\Sigma', t')\in \cD(\Xqm)$,
\[
\dd_{1, \theta}((\Sigma, t), (\Sigma', t'))=\vol(\theta; \Sigma, t)+\vol(\theta; \Sigma', t')-2\vol(\theta; \Sigma\cup\Sigma', t\wedge t'),
\]
where $t\wedge t'$ is the poinwise minimum on $\Sigma\wedge\Sigma'$ and is given by $t$ or $t'$ otherwise. For a big line bundle $L$, the map $\nn(L)$ therefore identifies $\cN^\qm(L)$ as the Hausdorff quotient of $(\cD(\Xqm), \dd_{1,\theta})$, and hence $(\cN^\qm(L), \de)$ also depends only on the numerical class of $L$, up to canonical homeomorphisms. Thus for $\theta\in \Bigc(X)$, even though norms are not defined for real classes, we can make sense of the metric space $(\cN^\qm(\theta), \de)$ as the Hausdorff quotient of $(\cD(\Xqm), \dd_{1,\theta})$. With this at hand, we now state the Calabi--Yau theorem for quasimonomial measures.

\bigskip\noindent\textbf{Theorem B.} 
For $\theta\in\Bigc(X)$, the Monge--Amp\`ere operator induces a bijection
\[
\MA_\theta\colon \cN^\qm(\theta)/\R\to \cM^\qm.
\]
Further, for each $\Sigma\subset\Xqm$ finite of size $d$, the Monge-Amp\`ere operator restricts to a homeomorphism
\[
\MA_\theta\colon (\cN^\Sigma(\theta)/\R, \qd_1)\to \sigma_\Sigma,
\]
where $\sigma_\Sigma\subset\R^d$ is the standard simplex equipped with the usual Euclidean topology, and $\qd_1$ is the quotient metric of $\de$.



\bigskip 

The key input in proving Theorem~B is again the differentiability of volumes in Theorem~A. Given $\mu=\sum_{i=1}^d\xi_i\delta_{v_i}\in \cM^\qm$, following the variational approach, a quasimonomial datum $(\Sigma, t)\in\cD(\Xqm)$ solves $\MA_\theta(\Sigma, t)=\mu$ iff it achieves the supremum in
\[
\vol^\vee(\theta; \Sigma, \xi)\=\sup_{t\in\R^d}\{\vol(\theta; \Sigma, t)-\langle \xi, t\rangle\},
\]
where $\xi\=(\xi_1, \dots, \xi_d)\in \sigma_\Sigma$ is the vector defined by $\mu$. We then get surjectivity of $\MA_\theta$ by restricting our variational problem to a compact subset of $\R^d$, which we get for free from the volume function. In contrast, in the usual approach of solving Monge--Amp\`ere equations through non-Archimedean pluripotential theory, such compactness in an infinite dimensional space of psh functions is often the hardest to establish, and relies on the \emph{envelope property}. 

The proof of injectivity consists of two steps. The first is the to analyze the support of the Monge--Amp\`ere measure $\MA_\theta(\Sigma,t)$. This is again in spirit close to the classical analysis of \emph{Rees valuations} for an integrally closed ideal. In our setting, we show that for a big line bundle $\theta=c_1(L)$ and $\chi\in \cN^\qm(L)$,
there is a unique minimal quasimonomial datum $(\rees(\chi), t)$ such that $\chi=\chi_{\rees(\chi),t}$. We further extend the set $\rees(\theta; \Sigma,t)$ to all of the big cone using $\dd_{1, \theta}$ pseudometrics. As a byproduct of the proof of Theorem~A, we get 
\[
\mathrm{supp}(\MA_\theta(\Sigma, t))=\rees(\theta; \Sigma, t).
\]
The second step exploits further the convexity properties of filtered Okounkov bodies, and for a fixed $\Sigma\subset\Xqm$, the uniqueness of solutions to $\MA_\theta(\Sigma, \cdot)=\mu$ is obtained through the Brunn--Minkowski inequality.

Theorem~B extends known results in the literature in various directions. When $\theta=c_1(L)$ is ample, this recovers~ \cite[Theorem~C]{bj:kstab1} for divisorial measures (see also~\cite[Theorem~7.13]{bj:kstab1} for a more general version, and~\cite[Proposition~8.6]{bfj:nama} for the non-trivially valued case).
For $\theta\in \Bigc(X)_\Q$, this extends the surjectivity result in~\cite{derreb:3} for divisorial date, 
where the Monge--Ampère operator is again defined using test configurations and therefore requires more machinery to deal with nonexistence of a test configuration with a given set of divisorial valuation as its precise Rees valuations; our proof is again simpler, and gives uniqueness conjectured in~\cite{derreb:3}.

\subsection*{Continuity of solutions with respect to big classes}
 Finally, an important advantage of working with quasimonomial data is that one can understand such a set independently of the choice of numerical classes. For a fixed finite set $\Sigma\subset\Xqm$, we can theorefore view the Monge--Amp\`ere measures as functions on the big cone. As a direct corollary of Theorem~B, we obtain
the following continuity result for solutions of non-Archimedean Monge--Ampère equations as the numerical class varies.

\bigskip\noindent\textbf{Theorem C.} 
    Let $(\mu_k)_k$ be a sequence of quasimonomial measures with support $\Sigma=\{v_i\}_{i=1}^d$ converging weakly to a quasimonomial measure $\mu$ with the same support. Assume $t_k\in \R^d$ solves $\MA_{\theta_k}(\Sigma, t_k)=\mu_k$, with $\theta_k\in \Bigc(X)$ converging to a class $\theta\in \Bigc(X)$. Then there exists a unique $t\in \R^d$ up to translation such that
    $$\MA_\theta(\Sigma,t)=\mu,\quad \qd_{1,\theta}((\Sigma,t_k),(\Sigma,t))\to 0,$$
where $\qd_{1, \theta}$ is the quotient pseudometric of $\dd_{1,\theta}$.

\bigskip 

We also remark that, due to the birational nature of constructions involving big classes, and birational invariance of $X^\qm$, most of our results in fact make sense when phrased in terms of \textit{big Cartier b-divisor classes} on the Riemann--Zariski space of $X$.
In particular, one could phrase Theorem C as obtaining continuity of solutions across varying birational models of $X$ as well. We do not pursue this approach here for clarity of writing.

\subsection*{Organization}
After recalling some background of norms and filtrations in~\S\ref{sec:prelim}, we define and study the space of valuation data/norms and their Rees valuations in~\S\ref{sec:valdata}, and extend the volume function to the big cone. In~\S\ref{sec:oko}, we study and extend Peng's Okounkov body construction with respect to a quasimonomial valuation in our setting. Finally we prove the main theorems in~\S\ref{sec:MA}, and compare our definition of the Monge--Amp\`ere operator to the existing ones in the literature.

\subsection*{Setup and notation}
\begin{itemize}
    \item We work over a trivially valued algebraically closed field $k$ of characteristic zero.
    \item $X$ is a normal projective variety over $k$ of dimension $n$.
    \item If $L$ is a big line bundle on $X$, we will write $R\=R(X, L)$ for the section ring, and $R_m\=H^0(X, mL)$ when no confusion arises.
    \item We write ${\bf 1}$ for the vector $(1,\dots, 1)\in \R^d$ when the dimension $d$ is clear from context.
    \item For a vector $\alpha\in \R^r$ and $p\in \R^n$ with $r\le n$, we write $\langle \alpha, p\rangle$ for the usual inner product $\langle \alpha, \pi(p)\rangle$ on $\R^r$, where $\pi\colon \R^n\to \R^r$ is the projection onto the first $r$ coordinates.  
    \item For a finite set $\Sigma$ of size $d$, $\sigma_\Sigma=\{x=(x_1, \dots, x_d)\in\R^d\mid \sum_{i=1}^dx_i=1\}$ is the standard simplex.
\end{itemize}

\subsection*{Acknowledgements.} We thank Ruadhaí Dervan and Lu Qi for discussions related to the contents of this article. We further thank Sébastien Boucksom and Mattias Jonsson for comments on a preliminary draft. The second author was supported in part by the Simons Collaboration Grant \emph{Moduli of varieties}.

\subsection*{AI usage.} The present project was started and mostly completed before the "AI era". Generative AI tools were used to provide the proof of Proposition~\ref{prop:d1uniqueness}, suggest the expression of $\FS$ in Corollary~\ref{cor:FSd1}, simplify the proof of Theorem~\ref{thm:diff} (which formerly used very technical convex geometric tools) and Theorem \ref{thm:realma}, and proofread. The paper was \textit{written} entirely by the authors.

\section{Preliminaries}\label{sec:prelim}
We start by recalling some basic facts about volumes of linear series and norms/filtrations, and Okounkov body constructions. The only new result in this section is Proposition~\ref{prop:fs} 
which should also be known to experts.
\subsection{Volume of line bundles and linear series}

A line bundle $L$ on $X$ is said to be \textit{big} if $\vol(L)>0$, where the volume is defined by
$$\vol(L)\=\lim_{m\to\infty}\frac{h^0(X,mL)}{m^n/n!}.$$
If $V\subset X$ is a subvariety of $X$ of dimension $p$, the \textit{restricted volume} of $L$ along $V$ is defined in~\cite{elmnp:restricted} as
$$\vol_{X|V}(L):=\lim_{m\to\infty}\frac{\dim \mathrm{Im}(H^0(X,mL)\rightarrow H^0(V,mL|_{V}))}{m^{p}/p!}.$$
In general, the volume function is continuously differentiable on the big cone, as shown in~\cite{bfj:teissier, lazmus}. If $X$ is smooth and $E\subset X$ is a prime divisor, then
$$\frac{d}{dt}\Big|_{t=0} \vol(L+tE) = n\vol_{X|E}(L).$$

More generally, we will consider graded linear series of a big line bundle $L$, i.e.\ graded subalgebras $V_\bullet\=(V_m)_m\subset R(X,L)$. The \textit{volume} of a graded linear series is defined as
$$\vol(V_\bullet):=\limsup_{m\to\infty}\frac{\dim(V_m)}{m^n/n!}.$$
Following \cite[Definition 1.1]{bouchen} we say $V_\bullet$ \textit{contains an ample series} if $V_m\neq 0$ for $m\gg 0$, and there exists $A$ an ample $\bbq$-Cartier $\bbq$-divisor with $L=A+E$, $E$ effective, such that $H^0(X,mA)\subseteq V_m (\subseteq H^0(X,mL))$ for all sufficiently divisible $m$. In particular, a graded linear series containing an ample series has positive volume.

\subsection{Berkovich analytification and valuations}
By a \emph{valuation} on $X$ we mean a real valuation $v\colon k(X)^\times\to \R$ extending the trivial valuation on $k$. We write $\Xval$ for the space of valuations on $X$, equipped with the topology of pointwise convergence on $k(X)^\times$. The space $\Xval$ admits a natural compactification by the \emph{Berkovich analytification} $\Xna$ of $X$. As a set, $\Xna\=\bigsqcup_V V^{\mathrm{val}}$, with $V$ ranging through all subvarieties $V\subseteq X$, and $\Xna$ is a compact Hausdorff topological space.

Given a valuation $v\in\Xval$, by the valuative criterion of properness, $v$ admits a \emph{center} $c_X(v)$ on $X$, as the unique scheme point $\xi$ such that $v\ge0$ on the local ring $\cO_{X,\xi}$ and $v>0$ on its maximal ideal. For $0\ne s\in H^0(X, L)$ a nonzero section of a line bundle $L$, $v(s)$ is the value of $v$ on the local function corresponding to $s$ under a trivialization of $L$ at $c_X(v)$.

As in~\cite{BKMS, blumjonsson}, the \textit{maximal vanishing order} of $v$ with respect to a big line bundle $L$ is
$$T(L;v):=\sup_m \sup\{m\mi v(s)\mid s\in H^0(X,mL)\setminus\{0\}\}\in (0, \infty].$$
More generally, we can define the maximal vanishing order $T(V_\bullet;v)$ for a graded linear series  $V_\bullet\subset R(X,L)$, for $L$ big, in the same way.
A valuation $v$ is \textit{linearly bounded} if $T(L;v)<\infty$ for some (or equivalently any) big line bundle $L$. We denote by $\Xlin\subset \Xval$ the set of linearly bounded.

The space $\Xlin$ contains the space $\Xdiv$ of divisorial valuations and more generally the space $\Xqm$ of \emph{quasimonomial valuations}. A valuation $v$ is \emph{quasimonomial} if it satisfies the Abhyankar equality
$$\mathrm{tr.deg}(v)+\mathrm{rat.rk}(v)=n,$$
where $\mathrm{tr.deg}(v)\=\mathrm{tr.deg}(k(v)/k)$ for $k(v)$ the residue fields of the valuation ring $\cO_v$, and $\mathrm{rat.rk}(v)\=\dim_\Q(v(K(X)^\times)\otimes_\Z \Q)$ is the \emph{rational rank} of $v$. A quasimonomial valuation is \emph{divisorial} if it has rational rank $\mathrm{rat.rk}(v)=1$. By convention, we also include the trivial valuation $v_{\triv}$ in $\Xdiv$, and we have inclusions of dense subsets
\[
\Xdiv\subset \Xqm\subset \Xlin\subset \Xna.
\]
By a theorem of Zariski~\cite[Lemma~2.45]{KM98}, a divisorial valuation is of the form $c \ord_E$ for $c\in\bbr_{>0}$, and $E$ a prime divisor in a normal birational model $\pi\colon Y\to X$. If $Y$ is further taken to be smooth, then
$$T(L;\ord_E)=\sup\{t>0\mid \pi^*L-tE\text{ is big}\}$$
coincides with the \textit{pseudoeffective threshold}.

More generally, as in~\cite{ELS, jonmus}, quasimonomial valuations can be characterized geometrically as follows.
Given a quasimonomial valuation $v$ of rational rank $r$, there is a log smooth model $(Y, D=\sum_j D_j)\to X$ and a stratum $\bigcap_{j=1}^rD_j\subset  Y$, with algebraic coordinates $(y_1,\dots,y_r)$ at the generic point $\eta$ of a connected component of $\bigcap_jD_j$, as well as $\alpha\in \bbr^r_{>0}$ consisting of $\Q$-linearly independent coordinates, such that for any $f\in \cO_{Y,\eta}$, if $f=\sum_\beta a_\beta y^\beta$ in $\widehat{\cO_{Y,\eta}}$, we have $v(f)=\min_{a_\beta\neq 0}\langle \alpha,\beta\rangle$. This property determines $v$. 
Given such a log smooth pair $(Y,D)$, we denote by $\QM(Y,D)$ the set of quasimonomial valuations obtained by strata of $(Y,D)$. Finally, we will often omit the log smooth pair in the notation and simply denote by $v=v_\alpha$ with $\alpha$ its weight.

\subsection{Norms and filtrations}
Fix a big line bundle $L$ and we write 
\[
R(X, L)\=\bigoplus_m R_m\=\bigoplus_mH^0(X, mL)
\]
for its section ring. By the previous discussion, any valuation $v\in \Xval$ defines a $k^\times$-invariant valuation on $R$. More generally, we consider \emph{norms} on $R$.
\begin{definition} A \emph{norm} on $R$ is a map $\chi: R\to \R\cup\{\infty\}$ with the following properties:
\begin{itemize}
    \item non-Archimedean: $\chi(f+g)\ge \min\{\chi(f), \chi(g)\}$ for all $f, g\in R$ and $\chi(cf)=\chi(f)$ for any $c\in k^\times, f\in R$;
    \item super-additive: $\chi(fg)\ge \chi(f)+\chi(g)$ for all $f, g\in R$;
    \item $k^\times$-invariant: $\chi(\sum_mf_m)=\min_m\chi(f_m)$, for $f=\sum_m f_m, f_m\in R_m$;
    \item $\chi(f)=\infty$ iff $f=0$.
\end{itemize}
A norm is further said to be \emph{linearly bounded} if there exists $C>0$ such that $|\chi|\le Cm$ on $R_m\setminus\{0\}$ for all $m\ge 1$.
\end{definition} 
The norm defined by a valuation $v$ is additive in the sense that $v(fg)=v(f)+v(g)$ for $f, g\in R$. Given any linearly bounded norm $\chi$, $c\in \R$ and $a>0$, we have that $a\chi$ is a norm, so is $\chi+c$, defined by
\[
(\chi+c)(f_m)\= \chi(f_m)+cm
\]
for any $f_m\in R_m$. If $\chi, \chi'$ are two norms, then their pointwise minimum $\chi\wedge\chi'$, given by 
\[
(\chi\wedge\chi')(f)=\min\{\chi(f), \chi'(f)\}
\]
for any $f\in R$ is again a norm.

The notion of maximal vanishing order generalizes to norms, by setting for each norm $\chi$, 
\[
T(\chi)\=\sup_m \sup\{m^{-1}\chi(s)\mid s\in R_m\setminus\{0\}\}.
\]
As observed in~\cite{bhj:duistermaat}, norms on $R$ are in one-to-one correspondence with filtrations on $R$.
\begin{definition}
    A \emph{filtration} $\cF$ on $R$ is given by a family of vector subspaces $\cF^\lambda R_m\subseteq R_m$ for $m\in \Z_{\ge0}, \lambda\in \R$, such that 
    \begin{itemize}
        \item $\cF^\lambda R_m\subseteq \cF^{\lambda'} R_m$ for any $\lambda\ge \lambda'$;
        \item $\cF^\lambda R_m=\bigcap_{\lambda'<\lambda} \cF^{\lambda'}R_m$;
        \item $\cF^\lambda R_m=R_m$ for $\lambda\ll 0$ and $\cF^\lambda R_m=0$ for $\lambda\gg 0$;
        \item $\cF^\lambda R_m\cdot \cF^{\lambda'}R_{m'}\subseteq \cF^{\lambda+\lambda'}R_{m+m'}$ for all $\lambda, \lambda'\in \R$, $m, m'\in \Z_{\ge0}$.
    \end{itemize}
\end{definition}
The correspondence between norms and filtrations on $R$ is obtained via
\[
\cF^\lambda_\chi R_m\=\{s\in R_m\mid \chi(s)\ge \lambda\}, \quad \chi_\cF(s)\=\sup\{\lambda\in \R\mid s\in \cF^\lambda R_m\}, s\in R_m,
\]

Given any norm $\chi$, or equivalently filtration $\cF_\chi$ on $R$ and $\lambda\in \R$, the associated graded linear series is defined as
\begin{equation}
V^\lambda_\bullet(\chi)\=V^\lambda_\bullet(\cF_\chi)\=\bigoplus_m V_m^\lambda\=\bigoplus_m\cF^{m\lambda}_\chi R_m.
\end{equation}
As in~\cite{bouchen, blumjonsson}, a filtration $\cF$ is \emph{linearly bounded} if the associated norm is linearly bounded, in particular,
\[
T(\cF)\=\lim_{m\to\infty}\sup\{\lambda\in \R\mid V_m^\lambda\ne 0\}=\sup\{\lambda\in\R \mid \vol(V_\bullet^\lambda(\cF))>0\}<\infty.
\]
and 
\[
\lambda_{\min}\=\inf\{\lambda\in\R\mid \vol(V^\lambda_\bullet(\chi))<\vol(L)\}>-\infty.
\]
It is easy to see that $T(\cF_\chi)=T(\chi)$ and so any $v\in\Xlin$ defines a linearly bounded filtration.
By~\cite[Lemma~1.6]{bouchen}, $V^\lambda_\bullet(\cF)$ contains an ample series for all $\lambda<T(\cF)$.



The \emph{volume} of a linearly bounded norm $\chi$ on $R$ (or equivalently, a linearly bounded filtration) is given by
\begin{equation}\label{eq:volint1}
\vol(\chi)\=c+\frac{1}{\vol(L)}\int_{c}^{T(\chi)} \vol(V^\lambda_\bullet(\chi))\, d\lambda,
\end{equation}
for any $c\le \lambda_{\min}$, see~\cite[Theorem~3.6]{bj:kstab1}.
Similar to the asymptotic characterization of $T(L; \chi)$, $\vol(\chi)$ can also be computed asymptotically, thanks to~\cite[Lemma~2.9]{blumjonsson}. For each $m\ge 1$, set
\[
\vol_m(\chi)=\frac{1}{mN_m}\sum_{i=1}^{N_m}\chi(s_i),
\]
for an orthogonal basis $(s_i)$ on $R_m$ for $\chi$. Then
\begin{equation}\label{eq:volasy}
    \vol(\chi)=\lim_{m\to\infty}\vol_m(\chi).
\end{equation}

Following work of Darvas~\cite{dar:mabuchi} in the complex analytic setting, a pseudodistance $\de$ was introduced on the space of norms on an ample line bundle~\cite{boueri}\cite[Section~3.4]{bj:kstab1} (which readily extends to the big case), by
\[
\de(\chi_1, \chi_2)=\vol(\chi_1)+\vol(\chi_2)-2\vol(\chi_1\wedge \chi_2),
\]
where $(\chi_1\wedge\chi_2)(s)=\min(\chi_1(s),\chi_2(s))$ denotes the pointwise minimum. It was shown in \cite[Corollary 6.24]{bj:kstab1} that if $L$ is ample, $\de$ restricts to a metric on maximal norms, and in particular $\cN^\div(L)$. We will later extend this result to the case of valuative norms on big line bundles.

\subsection{The Fubini--Study operator}
Following~\cite{bj:kstab1}, we also have the \emph{Fubini--Study operator} sending each linearly bounded norm $\chi$ on $R$ to a lower-semicontinuous function $\FS(\chi)\colon\Xna\to \R\cup \{-\infty\}$, by setting
\[
\FS(\chi)(v)\= \sup_m\frac 1m \sup\{\chi(s)-v(s)\mid s\in R_m\setminus\{0\}\}=\lim_m\frac 1m \sup\{\chi(s)-v(s)\mid s\in R_m\setminus\{0\}\},
\]
for each $v\in \Xna$. That the limit exists is a simple consequence of Fekete's lemma.
The following result should be well-known to experts, and the proof follows the same strategy that has already appeared multiple times in the literature in various settings, see~\cite{LX20, BLQ, jwu}.
\begin{proposition}\label{prop:fs}
    Let $\chi, \chi'$ be linearly bounded norms on $R=R(X, L)$ for $L$ big, such that $\chi\le \chi'$ and $\vol(\chi)=\vol(\chi')$. Then $\FS(\chi)=\FS(\chi')$ on $\Xlin$.
\end{proposition}
\begin{proof}
    We argue by contradiction, and suppose there is $v\in \Xlin$ such that $\FS(\chi)(v)<\FS(\chi')(v)$. Up to translating both $\chi, \chi'$, we may assume $\FS(\chi)(v)=0$, that is $\chi(s)\le v(s)$ for all $s\in R$ and that there exists a nonzero section $s\in R_m$ such that 
    \[
    0<\chi'(s)-v(s)\eqqcolon l.
    \]
    For any $p,d\in \Z_{>0}$, consider a map of the form
    \[
    \phi\colon R_d\to \cF^{\lambda(d+mp)}_{\chi'}R_{d+mp}/\cF^{\lambda(d+mp)}_{\chi}R_{d+mp},\quad  t\mapsto ts^p,
    \]
    for $\lambda$ to be determined. We aim to choose $\lambda$ so that $\phi$ is injective. To this end, note that since $\chi'$ is linearly bounded, there exists $C>0$ so that $\chi'(t)\ge -Cd$ for any $t\in R_d$, and we have
    \[
    \chi(ts^p)\le v(ts^p)=v(t)+pv(s)\le dT(v)+pv(s),
    \]
    and
    \[
    \chi'(ts^p)\ge \chi'(t)+p\chi'(s)=\chi'(t)+pl+pv(s) \ge -Cd+pl+pv(s).
    \]
    Thus for any $p\gg 0$ such that $pl>(T(v)+C)d$, and 
    \[
    \lambda(d+mp)\in (dT(v)+pv(s),-Cd+pl+pv(s)],
    \]
    $\phi$ is injective, and so 
    \[
    \dim \cF^{\lambda(d+mp)}_{\chi'}R_{d+mp}-\dim \cF^{\lambda(d+mp)}_{\chi}R_{d+mp}\ge \dim R_d,
    \]
    which then implies
    \[
    \vol(V^\lambda_\bullet(\chi'))-\vol(V^\lambda_\bullet(\chi))\ge \lim_{d\to \infty}\frac{d^n}{(d+mp)^n}\vol(L).
    \]
    Thus choosing $p=\lceil\frac{(T(v)+C)d}{l}\rceil+2$, we have
    \[
    \vol(\chi')-\vol(\chi)\ge \frac{l \lim_{d\to \infty}\frac{d^n}{(d+mp)^n}\vol(L)}{\vol(L)}>0,
    \]
    contradiction.
\end{proof}

\subsection{Okounkov bodies and concave transforms}

Okounkov bodies~\cite{okounkov,kkh,lazmus} are convex bodies capturing the volume of (a graded linear series of) a big line bundle, on which one can further define functions that capture the volume of linearly bounded filtrations~\cite{bouchen,wn:transforming}. Such Okounkov bodies are constructed as the convex hull of the rescaled semigroup of a $\mathbb{Z}^n$-valued valuation, for example coming from an admissible flag of $X$, see~\cite[Section 1.1]{lazmus}.

More precisely, given a graded linear series $W_\bullet\subset R(X,L)$ containing an ample series, and a valuation $\nu:W_\bullet\to\mathbb{Z}^n$ of maximal rational rank, let $\Gamma(W_\bullet)=\bigcup_m \{m\}\times \Gamma_m(W_\bullet)$ be the semigroup with 
$$\Gamma_m(W_\bullet)\=\nu(W_m\setminus\{0\}).$$
The Okounkov body $\Delta(W_\bullet)$ is the closed convex hull of the set of points $\bigcup_mm\mi \Gamma_m(W_\bullet)$, which is a convex body in $\bbr^n$. Then by~\cite{lazmus, bou:oko} one has 
$$\vol(W_\bullet)=n!\vol_{\bbr^n}(\Delta(W_\bullet)).$$
Note that this construction implicitly depends on the choice of such a valuation. Throughout, if $W_\bullet=R(X,L)$, we simply replace $W_\bullet$ with $L$ in our notation, so that e.g.\ $\Delta(L):=\Delta(R(X,L))$.

Now suppose we are given a norm $\chi$ on $R(X,L)$, recall that it induces for each $\lambda$ a graded linear series $V^\lambda_\bullet(\chi)\subset R(X,L)$ which contains an ample series for all $\lambda<T(\chi)$. Thus, for the rank $n$ valuation $\nu$ as above, we get the associated Okounkov body $\Delta(V^\lambda_\bullet(\chi))\subseteq \Delta(L)$ for each $\lambda<T(\chi).$ The \textit{concave transform} is
\begin{align*} 
G[\chi]\colon \Delta(L)&\to \bbr, \quad p\mapsto \sup\{\lambda\in \R\mid p\in\Delta(V^\lambda_\bullet(\chi))\}.
\end{align*}
By~\cite[Theorem 1.11]{bouchen} (see also~\cite{blumjonsson}), we have the following alternative formulas for $\vol(\chi)$:
\begin{equation}\label{eq:volint}
\vol(\chi)=\frac{1}{\vol(\Delta(L))}\int_{\Delta(L)}G[\chi](p)\,dp = \lambda_{\min}+\frac{1}{\vol(\Delta(L))}\int_{\lambda_{\min}}^\infty \vol(\Delta(V^\lambda_\bullet(\chi)))\,d\lambda.
\end{equation}
The concave transform can also be described as a \textit{concave envelope} of a function on the semigroup $\Gamma(L)$. For each $m$, and $p\in \Gamma_m(L)$ let
$$G_m[\chi](p):=\sup\{\chi(s)\mid\nu(s)=p,\,s\in R_m\setminus\{0\}\}.$$
Then we have that $G[\chi]$ is the concave envelope of the function $(m,p)\mapsto G_m[\chi](p)$ on the semigroup $\Gamma(L)$, in the sense of \cite[Definition 3.7]{wn:transforming}. In particular, by \cite[Theorem 3.8]{wn:transforming}, given $p\in \Delta(L)^\circ$ and a sequence $i\mapsto m_i\to\infty, p_i\in \Gamma_{m_i}(L)$ with $m_i\mi p_i\to p$, then
$$G[\chi](p)=\lim_{i\to\infty}m_i\mi G_{m_i}[\chi](p_i).$$

\section{The space of valuation data and valuative norms}\label{sec:valdata}
In this section, we define and study the space $\cD(\Xlin)$ of valuative data and the space of valuative norms $\cN^\lin(\theta)$ for a big class $\theta\in \Bigc(X)$. We begin by studying valuative norms on a big line bundle and identify a minimal set $\rees(\chi)$ of Rees valuations for each valuative norm. We then show that the volume function for valuative data is numerical in nature, and hence extends by continuity to all of $\Bigc(X)$. Through Darvas metrics, $\cN^\lin(\theta)$ is identified with the Hausdorff quotient of $\cD(\Xlin)$ for each $\theta$.

\subsection{Valuative data and valuative norms}
Let $L$ be a big line bundle on $X$, and write $R=\bigoplus R_m$ as before for the section ring of $L$. We first define the space of \emph{valuative data} and \emph{valuative norms}. When $L$ is ample, these norms and more generally maximal (or saturated) norms were studied in detail in~\cite{bj:kstab1} (see also~\cite{BLQ, jwu} for other settings). In this paper, however, we are only interested in norms defined by finitely many valuations on $X$.

\begin{definition}
    A \emph{valuation datum} is a pair $(\Sigma, t)$, where
\begin{itemize}
    \item $\Sigma=\{v_i\}_{i=1}^d\subset \Xlin$ is a finite set of valuations of linear growth on $X$;
    \item $t=(t_v)_{v\in\Sigma}\in\bbr^{d}$ is a vector.
\end{itemize}
A valuation datum is \emph{divisorial} (resp.\ \emph{quasimonomial}) if $\Sigma\subset X^\div$ (resp.\ $\Sigma\subset X^\qm$).
Denote by $\cD(\Xlin)$ the space of valuation data. We will also write $\cD(\Xdiv)$ (resp.\ $\cD(\Xqm)$) for the space of \emph{divisorial data} (resp. \emph{quasimonomial data}). For fixed $\Sigma\subset\Xlin$ of size $d$, the space $\cD(\Sigma)$ can be identified with $\R^d$.
\end{definition}


For each big line bundle $L$ on $X$, a valuation datum $(\Sigma, t)\in \cD(\Xlin)$ defines a \textit{valuative norm} $\chi_{\Sigma, t}\=\min_{v\in \Sigma}\{v+t_v\}$ on $L$ by
\[
\chi_{\Sigma,t}(s)=\min_{v\in \Sigma}\{v(s)+mt_v\},
\]
for any $s\in R_m$.
The filtration corresponding to a valuative norm defined by a valuation datum $(\Sigma,t)$ is then
$$\cF^\lambda_{\Sigma,t}R_m=\{s\in R_m \mid v(s)+mt_v\geq\lambda, \, \forall v\in\Sigma\}.$$

We say that $\chi_{\Sigma,t}$ is \textit{supported on $\Sigma$}, and denote by $\cN^\Sigma(L)$ the set of valuative norms supported on $\Sigma$. In a similar manner, we write $\cN^\lin(L), \cN^\div(L)$ and $\cN^\qm(L)$ for the space of valuative norms, divisorial norms and quasimonomial norms respectively. 
What we have described above is a natural map 
\[
\nn(L)\colon \cD(\Xlin)\to \cN^\lin(L),\quad (\Sigma, t)\mapsto \chi_{\Sigma, t},
\]
for any given big line bundle. This map is further $\R$-equivariant in the following sense. Given any $c\in \R$ and $(\Sigma, t)\in \cD(\Xlin)$ with $\Sigma$ having size $d$, we have
\[
\nn(L)(\Sigma, t+c{\bf 1})=\nn(L)(\Sigma, t)+c.
\]

In general, $\nn(L)$ is far from being injective, but we will see later that $\nn(L)$ admits a canonical minimal section. Indeed, given any (finite) subsets $\Sigma\subseteq \Sigma'\subseteq\Xlin$, it is easy to see that 
\[
T(\chi_{\Sigma, t})\le \inf_{v\in\Sigma} \{T(L; v)+t_v\},
\]
and so for any $w\in\Sigma'\setminus\Sigma$ and $t_w\ge \max_{v\in \Sigma}T(L; v)+\min_{v\in \Sigma} t_v$, we have
\[
\chi_{\Sigma, t}=\min\{\chi_{\Sigma, t}, w+t_w\}. 
\]
This gives an inclusion
\begin{equation}\label{eq:embedding}
\iota\colon\cN^\Sigma(L)\hookrightarrow\cN^{\Sigma'}(L), \quad \chi_{\Sigma, t}\mapsto \chi_{\Sigma', (t, c{\bf 1})},
\end{equation}
by choosing (any) $c=c(\Sigma, t)\gg 0$ for each $(\Sigma, t)$.
For later use, we also record the following convenient characterization of the minimum of two valuative norms, which follows immediately by unwinding definitions.
\begin{lemma}\label{lem:lattice}
    Let $\Sigma, \Sigma'\subset \Xlin$ be finite subsets of size $d$ and $d'$. Then for any $t\in \R^d$ and $t'\in \R^{d'}$, we have
    \[
    \chi_{\Sigma, t}\wedge \chi_{\Sigma', t'}= \chi_{\Sigma\cup \Sigma', t\wedge t'}
    \]
    where $(t\wedge t')_v=t_v\wedge t'_v$ if $v\in \Sigma\cap \Sigma'$, $(t\wedge t')_v=t_v$ if $v\in \Sigma\setminus\Sigma'$ and $(t\wedge t')_v=t'_v$ otherwise.    
    In particular, if $\chi=\chi_{\Sigma, t}=\chi_{\Sigma, t'}$ for $t,t'\in \R^d$, then
    $$\chi=\chi_{\Sigma,t\wedge t'}.$$
\end{lemma}

A crucial point of view taken throughout this paper is that valuation data give norms on each big line bundle on $X$, and thus we may understand invariants of valuative norms as line bundles vary in the same or different numerical classes. For a fixed big line bundle $L$, we may define invariants on $\cD(\Xlin)$ by pulling back the invariants on $\cN^\lin(L)$ via $\nn(L)$. More precisely, we will write
\begin{align*} 
    T(L;\Sigma,t)&\=T(\nn(L)(\Sigma, t));\\
    \vol(L;\Sigma,t)&\=\vol(\nn(L)(\Sigma, t));\\
    V^\lambda_\bullet(L;\Sigma,t)&\=V^\lambda_\bullet(\nn(L)(\Sigma, t));\\
    \FS(L; \Sigma, t)&\= \FS(\nn(L)(\Sigma, t))\\    
\end{align*}
More generally, since valuation data also give norms on graded linear series $W_\bullet\subseteq R(X,L)$, we will also occasionally write e.g.\ $T(W_\bullet;\Sigma,t)$ for $T(\chi_{\Sigma,t}(W_\bullet))$.

\begin{remark} 
    We recover usual invariants of valuations by taking $\Sigma$ to be a singleton. For example, $T(L;v)=T(L;\{v\},0)$. Likewise, $\vol(L;\{v\},0)$ is often written as $S_L(v)$ in the literature. When $\Sigma$ is more generally a finite subset of $X^\div$ associated with divisors $(E_i)_i$, $\vol(L;\Sigma,t)$ is also written as $S_L(\{E_i,t_i\})$ in \cite{derreb:3}.
\end{remark}

In particular, for each big line bundle $L$, we get a pseudometric $\dd_{1, L}$ on $\cD(\Xlin)$ by setting for any $(\Sigma, t), (\Sigma', t')\in \cD(\Xlin)$,
\[
\dd_{1,L}((\Sigma, t), (\Sigma', t'))\=\de(\chi_{\Sigma, t}, \chi_{\Sigma', t'}).
\]
More concretely, using Lemma~\ref{lem:lattice} we have
\begin{equation}\label{eq:d1valdata}
\dd_{1,L}((\Sigma,t),(\Sigma',t'))=\vol(L;\Sigma,t)+\vol(L;\Sigma',t')-2\vol(L;\Sigma\cup\Sigma',t\wedge t').
\end{equation}
The natural map $\nn(L)$ upgrades to a distance preserving map 
\[
\nn(L)\colon (\cD(\Xlin), \dd_{1,L})\to (\cN^\lin(L), \de)
\]
between pseudometric spaces.

\subsection{Rees valuations}
We now explain that there exists a unique minimal representation of a valuative norm $\chi$ supported on the set $\rees(\chi)$ of Rees valuations. This is in turn inspired by the theory of Rees valuations for integrally closed ideals.
\begin{definition}
    Let $\chi=\chi_{\Sigma, t}$ for a finite set $\Sigma\subset\Xlin$. We say that $\Sigma$ is \emph{irredundant} for $\chi$ if for any $v\in \Sigma$, there is some $s\in R_m$ such that $\chi(s)=v(s)+mt_v<\min_{w\in\Sigma\setminus\{v\}}\{w(s)+mt_w\}$.
\end{definition}
\begin{definition}
    Let $\chi\in\cN^\lin(L)$. The set of \emph{Rees valuations} of $\chi=\chi_{\Sigma, t}$ is defined as
    \[
    \rees(\chi)\=\rees(L; \Sigma, t) = \bigcap\{\Sigma'\subseteq \Sigma\mid \chi_{\Sigma', t}=\chi_{\Sigma, t} \}
    \]
    where we have used the same $t$ for the restriction of $t$ to $\Sigma'$.
\end{definition}

Recall that classically, for an ideal $\fa\subset\cO_X$, the set of Rees valuations of $\fa$ is the unique irredundant set $S$ of valuations satisfying $\overline{\fa^m}=\bigcap_{v\in S}\{f\mid v(f)\ge mv(\fa)\}$, see e.g.~\cite[Theorem~10.1.6]{HS}. Our next result shows that the set $\rees(L; \Sigma,t)$ indeed gives the minimal irredundant set for a valuative norm, and in particular only depends on the norm. As a first step, we first show that the irredundant set is unique for a single valuation.
\begin{lemma}\label{lem:valrepunique}
    Suppose $v\in \Xlin$ and $\Sigma\subset \Xlin$ is a finite subset such that $v=\inf_{w\in \Sigma}\{w+t_w\}$ 
    and $\Sigma$ is irredundant for $v$, then $\Sigma=\{v\}$.
\end{lemma}
\begin{proof}
    Note that if $v\in \Sigma$, then since $\Sigma$ is irredundant for $v$, we necessarily have $t_v=0$, from which it follows that $\Sigma=\{v\}$. Suppose now that $v\notin \Sigma$, we will derive a contradiction.

    We first observe that if $\Sigma=\{w\}$, then $v=w$ and $t_w=0$. 
    Indeed, if $L$ is ample, then for $m\gg 0$, we may find $s\in H^0(X, mL)$ such that $0=v(s)=w(s)+mt_w$, and thus $t_w\le 0$, and we have $t_w=0$ by symmetry. In particular, $v=w$. In general, pick $l\gg 0$ and write $L'\= lL=A+E$ for an ample integer divisor $A$ and effective divisor $E$. Then by assumption we have $v=w+kt_w$ on $R(X, L')$. On the other hand, for any $s\in H^0(X, mA)$, multiplication by $s_E$ gives $ss_E^m\in H^0(X, mL')$, and we get
    \[
    v(ss_E^m)=w(ss_E^m)+mlt_w,
    \]
    and therefore $v=w+lt_w+w(s_E)-v(s_E)$ on $R(X,A)$. By the previous argument one has $v=w$, and in particular $v(s_E)=w(s_E)$ and $t=0$.
    
    Suppose now there is some $w'\in \Sigma$ with $w'\neq w$. By irredundancy of $\Sigma$, we can find $s\in R_m, s'\in R_{m'}$ such that 
    \begin{align*}
    v(s)&=w(s)+mt_w<\inf_{u\in \Sigma\setminus\{w\}}\{u(s)+mt_u\}, \\
    v(s')&=w'(s')+m't_{w'}<\inf_{u\in \Sigma\setminus\{w'\}}\{u(s)+m't_u\}.
    \end{align*}
    Thus for $N\gg 0$, we have
    \begin{align*}
    w'(s^Ns')+(mN+m')t_{w'}&=N(w'(s)+mt_{w'})+w'(s')+m't_{w'}\\
    &> N(w(s)+mt_w)+w(s')+m' t_w\\
    &=w(s^Ns')+(mN+m')t_w.
    \end{align*}
    Further, for any $u\in \Sigma\setminus\{w,w'\}$, we have
    \[
    u(s^Ns')+(mN+m')t_u=N(u(s)+mt_u)+(u(s')+m't_u))> Nv(s)+v(s')=v(s^Ns').
    \]
    It follows that 
    \[
    v(s^Ns')=\min\{w(s^Ns')+(mN+m')t_w,
     w'(s^Ns')+(mN+m')t_{w'}\}=w(s^Ns')+(mN+m')t_w,
    \]
    and this is a contradiction since $v(s')<w(s')+m't_{w}$.
\end{proof}

\begin{proposition}\label{prop:reesnorm}
    Let $\chi\in\cN^\lin(L)$ be a valuative norm.
    Then
    \begin{enumerate}[(i), font=\normalfont]
        \item if $\chi\in\cN^\Sigma(L)$, then $\chi=\chi_{\Sigma, t}$ with $t_v=\FS(L;\Sigma,t)(v)$. Further $t$ is minimal in the sense that, for any $t'$ such that $\chi=\chi_{\Sigma, t'}$, we have $t_v\le t_v'$ for all $v\in \Sigma.$  
        \item if $\Sigma$ is irredundant for $\chi=\chi_{\Sigma, t}$, then $t$ is minimal;
        \item $\rees(L; \Sigma, t)$ is the unique minimal irredundant set for $\chi=\chi_{\Sigma,t}$, that is, 
        for any $(\Sigma', t')\in\cD(\Xlin)$ with $\chi=\chi_{\Sigma',t'}$, 
         we have $\rees(\chi)\subseteq \Sigma'$ and $t_v\le t_v'$ for any $v\in \rees(\chi)$.
    \end{enumerate}
\end{proposition}
\begin{proof}
    The first assertion follows from~\cite[Proposition~6.5(i)]{bj:kstab1}, which is stated for $L$ ample; for completeness, we recall the argument here. 
    Assume $\chi=\chi_{\Sigma, t'}$. For each $s\in H^0(X, mL)$ and $v\in \Sigma$, we have
    $$\chi(s)=\min_{w\in\Sigma}(w(s)+mt_w')\leq v(s)+mt_v',$$
    that is, $m\mi(\chi(s)-v(s))\leq t_v'$ for all $s\in R_m$ and for all $m$, whence $t_v\leq t_v'$. This implies $\chi_{\Sigma,t}\leq \chi_{\Sigma,t'}$. Conversely, for $v$ and $s$ as above, we always have
    $$v(s)+m t_v\geq v(s)+m(m\mi(\chi(s)-v(s)))=\chi(s)$$
    which, after taking minimum over $\Sigma$, gives $\chi_{\Sigma,t}\geq\chi_{\Sigma,t'}$, hence equality.

    For~(ii), assume $\Sigma$ is irredundant, and assume to the contrary that $t$ is not minimal, and $\chi=\chi_{\Sigma, t'}$ for some $t'\neq t$. By Lemma~\ref{lem:lattice}, after replacing $t'$ by $\min\{t_v, t_v'\}$, we may assume $t'_v<t_v$ for some $v\in \Sigma$. By irredundancy, there exists $s\in R_m$ such that $\chi(s)=v(s)+mt_v$. But then
    \[
    \chi(s)\le v(s)+mt_v'<v(s)+mt_v=\chi(s),
    \]
    contradiction.

    Finally we prove~(iii). We first show that $v\in \Sigma$ is redundant for $\chi$ ,  i.e., $\chi_{\Sigma\setminus\{v\}, t\setminus\{t_v\}}=\chi_{\Sigma, t}$ iff $v\not\in\rees(L; \Sigma, t)$. Indeed, if $v$ is redundant, then $\rees(L; \Sigma, t)\subseteq \Sigma\setminus\{v\}$ by definition; conversely, given $v\notin\rees(L; \Sigma, t)$, there exists $\Sigma'\subseteq \Sigma$ with $\chi_{\Sigma', t}=\chi_{\Sigma, t}$ and $v\notin \Sigma'$. This implies $\Sigma'\subseteq \Sigma\setminus\{v\}$ and so 
    \[
    \chi_{\Sigma\setminus\{v\}, t\setminus\{t_v\}}\le \chi_{\Sigma', t}=\chi_{\Sigma,t}.
    \]
    This implies $\chi=\chi_{\rees(L; \Sigma, t), t}$ and that $\rees(L; \Sigma,t)$ is irredundant for $\chi$, since 
    \[
    \chi=\min_{v\in\Sigma\setminus\rees(L; \Sigma, t)}\{\chi_{\Sigma\setminus\{v\}, t\setminus\{t_v\}}\}\ge \chi_{\rees(L; \Sigma, t), t}\ge\chi.
    \]
    It remains to show that if $\chi=\chi_{\Sigma', t'}$, then $\rees(L; \Sigma, t)=\rees(L; \Sigma', t')$ and $t_v= t_v'$ for $v\in \rees(L; \Sigma, t)$. After replacing $\Sigma$ and $\Sigma'$ by their set of Rees valuations, we will assume that both $\Sigma$ and $\Sigma'$ are irredundant for $\chi$, and we aim to show $\Sigma\subseteq\Sigma'$.
    Pick $v\in \Sigma$, and $s_0\in R_m$ such that 
    \[
    \chi(s_0)=v(s_0)+mt_v<\min_{v'\in \Sigma\setminus\{v\}}\{v'(s)+mt_{v'}\}.
    \]
    Denote by $\Sigma''\=\{w\in \Sigma'\mid \chi(s_0)=w(s_0)+mt'_w\}$. Now let $s\in R_l\setminus\{0\}$ be any section, and $w\in \Sigma'\setminus\Sigma''$. First note that there exists $N_0\=N_0(v, s_0, s, \Sigma)$ such that for any $v'\in \Sigma\setminus\{v\}$ and $N\ge N_0$, we have
    \begin{align*}
    &v(ss_0^N)+(mN+l)t_v=N(v(s_0)+mt_v)+v(s)+lt_v\\
    &<N(v'(s_0)+mt_v)+v'(s)+kt_{v'}=v'(ss_0^N)+(mN+l)t_{v'},
    \end{align*}
    and in particular $\chi(s's_0^N)=v(ss_0^N)+(mN+l)t_v$. Similarly, we may find $N_1=N_1(s_0, s, v, w)>N_0$ such that 
    \begin{align*}
    w(s_0^Ns)+(Nm+k)t_w'&=N(w(s_0)+t_w)+w(s)+lt_w'\\
    &>N(v(s_0)+mt_v)+v(s)+lt_v=v(s_0^Ns')+(Nm+l)t_v,
    \end{align*}
    for all $N\ge N_1$. Taking $N_2\=\max_{w\in \Sigma'\setminus\Sigma''}N_1(s_0, s, v, w)$, we get that for $N\ge N_2(s_0, s, v, \Sigma)$, 
    \[
    \chi(s_0^Ns)=v(ss_0^N)+(mN+l)t_v=\min_{w\in \Sigma''}\{w(s_0^Ns)+(Nm+l)t_w'\}.
    \]
    By construction of $\Sigma''$, this yields
    \[
    v(s)+lt_v=\min_{w\in \Sigma''}\{w(s)+kt_w\}.
    \]
    Since $s$ is arbitrary, we get $v=\min_{w\in\Sigma''}\{w+t_w-t_v\}$, and by Lemma~\ref{lem:valrepunique} we have $\Sigma''=\{v\}$. This completes the proof of~(iii) in view of~(ii).

\end{proof}


As a consequence of Proposition~\ref{prop:reesnorm}, we get a canonical section of $\nn(L)$, by setting
\[
\rees\colon \cN^\lin(L)\to \cD(\Xlin), \chi\mapsto (\rees(\chi), \FS(\chi)|_{\rees(\chi)}).
\]
It also follows immediately from the definition that for any $t>0$, $c\in \R$, we have
\[
\rees(t\chi)=t\rees(\chi),\quad \rees(\chi+c)=\rees(\chi).
\]

\begin{proposition}\label{prop:reescvx}
    For any $\Sigma\subset\Xlin$ a finite subset of size $d$, the set 
    \[
    U(\Sigma)\=\{t\in \R^{d}\mid \rees(\chi_{\Sigma, t})=\Sigma\}\subset\R^d
    \]
    is open and convex. 
\end{proposition}
\begin{proof}
    Let $t\in U(\Sigma)$. 
    By assumption, for any $v\in \Sigma$, we can find $s\in R_m$ such that 
    \[
    v(s)+mt_v<\min_{w\in \Sigma\setminus\{v\}}\{w(s)+mt_w\}.
    \]
    This is an open condition on $t$, and thus still holds for $0<|t'-t|\ll1$. Since $\Sigma$ is finite, after taking finitely many intersections, there is an open neighborhood $U_t$ of $t$ so that $\Sigma$ is still irredundant for $\chi_{\Sigma, t'}$ for all $t'\in U_t$.

    Suppose now $t, t'\in U(\Sigma)$. We need to show $\lambda t+(1-\lambda) t'\in U(\Sigma)$ for all $\lambda\in [0,1]$. By openness, it suffices to assume $\lambda=\frac pq$ for positive integers $0<p<q$. To this end, pick $s\in R_m, s'\in R_{m'}$ such that 
    \[
    v(s)+mt_v<\min_{w\in \Sigma\setminus\{v\}}\{w(s)+mt_w\}, \quad
    v(s')+m't_v'<\min_{w\in \Sigma\setminus\{v\}}\{w(s')+m't_w'\}
    \]
    Then for $a\=pm', b\=(q-p)m$, we have
    \[
    v(s^as'^b)+amt_v+bm't_v'<w(s^as'^b)+amt_w+bm't_w'
    \]
    for any $w\in \Sigma\setminus\{v\}$. On the other hand, by our choice of $a,b$, we have
    \[
    amt_v+bm't_v'=(am+bm')\left(\frac{am}{am+bm'}t_v+\frac{bm'}{am+bm'}t_v'\right)
    =(am+bm')(\lambda t_v+(1-\lambda)t_v').
    \]
    This shows 
    \[
    v(s^as'^b)+(am+bm')(\lambda t_v+(1-\lambda)t_v')<\min_{w\in \Sigma\setminus\{v\}}\{w(s^as'^b)+(am+bm')(\lambda t_w+(1-\lambda)t_w')\},
    \]
    and so $\lambda t+(1-\lambda)t'\in U(\Sigma)$ as desired.
\end{proof}

\subsection{Darvas metrics and geodesics}

Turning to the metric structures, we first state the following analog of~\cite[Theorem~D]{bj:kstab1}.
\begin{theorem}\label{thm:d1metric}
    For any big line bundle $L$, the map  $\nn(L)\colon \cD(\Xlin)\to \cN^\lin(L)$ identifies $(\cN^\lin(L), \de)$ as the Hausdorff quotient of $(\cD(\Xlin), \dd_{1,L})$.
\end{theorem}
\begin{proof}
    It is clear from the definition that $\dd_{1, L}\equiv 0$ along fibers of $\nn(L)$. It remains to show that $\de$ is a metric on $\cN^\lin(L)$. For this, suppose $\chi, \chi'\in \cN^\lin(L)$ are such that $\de(\chi, \chi')=0$. Upon replacing one of them by $\chi\wedge \chi'$, we may assume $\chi\le \chi'$. Then we have $\vol(\chi)=\vol(\chi')$, and so $\FS(\chi)=\FS(\chi')$ on $\Xlin$ by Proposition~\ref{prop:fs}. If $\chi\in \cN^\Sigma(L)$ and $\chi'\in \cN^{\Sigma'}(L)$, then through the embedding described before, we have $\chi, \chi'\in \cN^{\Sigma\cup\Sigma'}(L)$.    
    It then follows from Proposition~\ref{prop:reesnorm}(i) that $\chi=\chi'$ as desired.
\end{proof}

Together with Proposition~\ref{prop:reesnorm}, this gives alternative characterizations of Rees valuations and the Fubini--Study operator on $\Xlin$.
\begin{corollary}\label{cor:FSd1}
    Given $(\Sigma, t)\in \cD(\Xlin)$ and a big line bundle $L$, we have
     \[
    \rees(L; \Sigma, t)=\bigcap\{\Sigma'\subseteq \Sigma\mid \vol(L; \Sigma', t)=\vol(L; \Sigma, t)\} = \bigcap\{\Sigma'\subseteq \Sigma\mid \dd_{1,L} ((\Sigma', t),(\Sigma, t))=0\},
    \]
    where we have used the same $t$ for the restriction of $t$ to $\Sigma'$, and     
    $$\FS(L;\Sigma,t)(v)=\inf\{c\in\bbr\mid \dd_{1,L}((\Sigma,t),(\Sigma\cup\{v\},t\cup \{c\}))=0\}$$
    on $\Xlin$, where in the above expression, if $v\in \Sigma$, the corresponding coordinate $t_v$ of $t$ is replaced with $t_v\wedge c$.
\end{corollary}
\begin{proof}
    The first statement follows immediately from Theorem~\ref{thm:d1metric}.
    Put $\chi=\chi_{\Sigma,t}$ and $\chi_c=\min(\chi,v+c)$. Then by the previous theorem, $\de(\chi, \chi_c)=0$ iff $\chi=\chi_c$. On the other hand, Proposition~\ref{prop:reesnorm}(i) yields
    \[
    \inf\{c\in \R\mid \chi=\chi_c\}=\FS(\chi)(v).
    \]
\end{proof}

As in \cite{reb:2,blxz}, given two norms $\chi_0,\chi_1$ on $R(X,L)$, with $L$ big, we can construct a geodesic segment $[0,1]\ni t\mapsto \chi_t$ for $\de$ by fixing in each degree a basis $(s_{m,i})_i$ of $R_m$ adapted to $\chi_0$ and $\chi_1$ and setting, for any $s=\sum_i a_{m,i} s_{m,i}\in R_m$,
$$\chi_t(s)=\min_{i,\,a_{m,i}\neq 0} \{(1-t)\chi_0(s_{m,i})+t\chi_1(s_{m,i})\}.$$
It follows from~\eqref{eq:volasy} (see also~\cite{reb:2,blxz}) that
\begin{equation}\label{eq:affine}
\vol(\chi_t)=(1-t)\vol(\chi_0)+t\vol(\chi_1).
\end{equation}

\begin{remark} 
A geodesic between two valuative norms supported on some $\Sigma\subset \Xlin$ does not necessarily remain in $\cN^\Sigma(L)$. In particular, it is not the case in general that the geodesic between $\chi_{\Sigma,t_0}$ and $\chi_{\Sigma,t_1}$ at time $t$ is given by $\chi_{\Sigma,(1-t)t_0+tt_1}$.
Let $v_1, v_2$ be monomial valuations with weights $(1,0), (0,1)$ on an affine patch in $\bbp^2$, and put $\chi_1=\min\{v_1+1, v_2\}, \chi_2=\min\{v_1, v_2+1\}$. Then $\chi_1=\chi_2$ on $R_2$ with the standard monomial basis being orthogonal. On the other hand, $\chi_3\=\min\{v_1+\frac 12, v_2+\frac12\}$ satisfies $\chi_3(z_0z_1)=2>\frac12\chi_1(z_0z_1)+\frac 12\chi_2(z_0z_1)=1$.
\end{remark}

\subsection{The volume with respect to big classes}

We fix in this section a finite set $\Sigma=\{v_i\}_{i=1}^d\subset \Xlin$ of valuations on $X$. For any big line bundle $L$  on $X$, the volume of valuative norms supported on $\Sigma$ induces a function
\[
\vol(L; \Sigma, \cdot)\colon \R^d\to \R,
\]
that is increasing in each variable, satisfies $\vol(L; \Sigma, t+c{\bf 1})=\vol(L; \Sigma, t)+c$, where ${\bf 1}=(1,\dots, 1)\in \R^d$, for any $t\in \R^d, c\in \R$. These two properties formally imply that the volume function is $1$-Lipschitz with respect to the $L^\infty$ norm on $\R^d$. Since the volume only depends on the norm $\chi_{\Sigma, t}$, it is the same for valuation data defining the same norm. In particular, $\vol(L; \Sigma, t)=\vol(L; \rees(L; \Sigma, t), t)$. Moreover, if $\min_{v_i\in \Sigma}t_i=0$, then for any $v_j\in \Sigma$ with $t_j\geq \max_{v_i\in\Sigma} T(L;v_i)$ and all $x\ge 0$, we have 
\begin{equation}\label{eq:volconst}
    \vol(L; \Sigma, t)=\vol(L; \Sigma, t+xe_j),
\end{equation}
where $e_j$ denotes the $j$-th unit vector in $\R^d$.

We now extend the volume function $L\mapsto \vol(L; \Sigma, t)$ to a continuous function on the big cone, following closely the classical theory of the volumes of divisors, see~\cite[\S2.2.C]{book:laz1}.

\begin{lemma}\label{lem:ratclass}
Let $r\in \Z_{>0}$. Then
\begin{enumerate}[(i), font=\normalfont]
    \item $\vol(r L;\Sigma,r t)=r\vol(L;\Sigma,t)$;
    \item $T(rL; \Sigma, r t)=r T(L; \Sigma, t)$;
    \item for fixed $(\Sigma, t)$, functions $\vol(L; \Sigma, t), T(L; \Sigma, t)$ only depend on the numerical class of $L$.
\end{enumerate}
\end{lemma}
\begin{proof}
    The first assertion follows from a direct calculation using~\cite[Lemma~2.9]{blumjonsson}. Indeed, for any $m\in \Z_{>0}$, put $N_m\=h^0(X, mL)$, and pick an orthogonal basis $\{s_i\}$ for $\chi_{\Sigma, rt}$, then we have
    \[
    \vol_m(rL; \Sigma, r t)=\frac{1}{mN_{mr}}\sum_{i=1}^{N_{mr}}\min_{v\in \Sigma}(v(s_i)+mrt_v) = r \vol_{mr}(L; \Sigma, t).
    \]
    Taking $m\to \infty$ yields the desired result. Similarly,~(ii) follows from the above asymptotics as well.
    
    The last statement will follow if it is true for $\min_v t_v=0$. In this case the same proof as in~\cite[Proposition~3.1]{BKMS} (see also~\cite[Proposition~4.1(i)]{lazmus}) yields the result, which we include here for the reader's convenience. Let $L'\=L+P$ for a numerically trivial line bundle $P$. By~\cite[Lemma~2.2.42]{book:laz1} we also fix a very ample line bundle $A$ such that $A+lP$ is very ample for any $l\in \Z$. Write $\chi\=\chi_{\Sigma, t}$ and fix a rank $n$ valuation $\nu$ of maximal rational rank giving rise to Okounkov bodies $\Delta(L)$ and $\Delta(L')$. We will show that for any $0\le \lambda<T(L; \Sigma, t)$, 
    \[
    \Delta(L; V^\lambda_\bullet(\chi))=\Delta(L'; V^\lambda_\bullet(\chi)).
    \]
    Since $L$ is big, there exists $k\gg1$ and a nonzero section $\sigma\in H^0(X, kL-A)$. We can write
    \[
    (m+k)L'=mL+(kL-A)+(A+(m+k)P).
    \]
    Since $A$ is very ample, there exists $0\ne \tau_m\in H^0(X, A+(m+k)P)$ for each $m\ge0$ such that $\tau_m$ does not vanish at the center of $\nu$, so that $\nu(\tau_m)=0$. Thus for any $s\in \cF^{\lambda m}_\chi H^0(X, mL)$, if $s'\=s\cdot \sigma\cdot \tau_m$, then $s'\in H^0(X, (m+k)L')$, $\chi(s')\ge \chi(s)\ge \lambda m$ and $\nu(s')=\nu(s)+\nu(\sigma)$. Thus for any $0<\varepsilon\ll1$ and $m\gg1$ such that $\lambda-\varepsilon<\lambda\frac{m}{m+k}<\lambda$, we have
    \[
    \nu(\cF^{\lambda m}_\chi H^0(X, mL))+\nu(\sigma)\subseteq \nu(\cF_\chi^{(\lambda-\varepsilon)(m+k)}H^0(X, (m+k)L').
    \]
    This shows 
    \[
    \Delta(L; V^\lambda_\bullet(\chi))\subseteq \bigcap_{0<\varepsilon\ll1}\Delta(L'; V^{\lambda-\varepsilon}_\bullet(\chi))=\Delta(L'; V^\lambda_{\bullet}(\chi)),
    \]
    and we conclude by symmetry.
\end{proof}

This allows us to define $\vol(\cdot; \Sigma, t)$ and $T(\cdot; \Sigma, t)$ on $\Bigc(X)_\Q$. For $L$ a big $\Q$-line bundle and $\theta=c_1(L)$, we set
\[
T(\theta; \Sigma, t)\=r^{-1}T(rL; \Sigma, r t), \quad \vol(\theta; \Sigma, t)\=r^{-1}\vol(rL; \Sigma, r t)
\]
for $r$ divisible enough. To define the volume function on all of $\Bigc(X)$, we need the following continuity results.
\begin{proposition}\label{prop:Tcts}
    The function $\theta\mapsto T(\theta;\Sigma, t)$ is concave on $\Bigc(X)_\Q$, and hence extends uniquely to a continuous concave function on $\Bigc(X)$.
\end{proposition}
\begin{proof}
    Fix $c\in[0,1]\cap \Q$, and $\theta, \theta'\in \Bigc(X)_\Q$.
    After possibly passing to a large multiple, we may assume both $(1-c)\theta$ and $c\theta'$ are represented by line bundles $(1-c)L$ and $cL'$. Then by definition we have
    \begin{align*}
    T((1-c)\theta+c\theta'; \Sigma, t)
    &\ge \sup_{D\in |L|_\Q, D'\in |L'|_\Q} \min_{v\in \Sigma}\{v((1-c)D+cD')+t_v\}\\
    &\ge \sup_{D\in |L|_\Q, D'\in |L'|_\Q}(1-c)\min_{v\in \Sigma}\{v(D)+t_v\}+c\min_{v\in \Sigma}\{v(D')+t_v\}\\
    &=(1-c)T(\theta; \Sigma, t)+cT(\theta'; \Sigma, t).
    \end{align*}
    The second claim follows from concavity.
\end{proof}
\begin{proposition}\label{prop:volcts}
    The function $\theta\mapsto \vol(\theta; \Sigma, t)$ on $\Bigc(X)_\Q$ extends uniquely to a continuous function on $\Bigc(X)$.
\end{proposition}
\begin{proof}
    Pick very ample divisors $A_1, \dots, A_r$ on $X$ whose classes form a basis for $\mathrm{NS}(X)_\Q$ and fix a norm $\|\cdot\|$ on $N^1(X)_\R$. We first assume $L=\sum_ia_iA_i, L'=\sum_ia_i'A_i$ with $a_i, a_i'\in \Z$ and $a_i\ge a_i'$. Put $B\=L-L'\ge0$, and pick $E_j\in |A_j|$ very general so that $s_{E_j}$ does not vanish at the center of $v\in \Sigma$. Now for each $m\in \Z_{>0}$ and $\lambda\in\R$, the multiplication map $H^0(X, mL')\to H^0(X, mL), s\mapsto s\prod_{j=1}^rs_{E_j}^{m(a_j-a_j')}$ induces an injection 
    \[
    \cF^{m\lambda}_{\Sigma, t}H^0(X, mL')\to \cF^{m\lambda}_{\Sigma, t}H^0(X, mL).
    \]
    The same argument as in the proof of~\cite[Theorem~2.2.44]{book:laz1} yields
    \[
    0\le\vol_L(V^{\lambda}_\bullet(\chi_{\Sigma, t}))-\vol_{L'}(V^{\lambda}_\bullet(\chi_{\Sigma, t}))\le \frac C2 \|L\|^{n-1}\|B\|,
    \]
    for some $C>0$. Similarly, if $L'=L+B-B'$ for $B, B'$ effective integral linear combinations of $A_j$, then 
    \[
    |\vol_L(V^{\lambda}_\bullet(\chi_{\Sigma, t}))-\vol_{L'}(V^{\lambda}_\bullet(\chi_{\Sigma, t}))|\le C\max\{\|L\|, \|L'\|\}^{n-1}\|L-L'\|.
    \]
    After possibly translating $t$ so that we can take $c=0$ in~\eqref{eq:volint1} and integrating over $\lambda$, we get
    \begin{align*}
    &|\vol(L)\vol(L; \Sigma, t)-\vol(L')\vol(L'; \Sigma, t)|\\
    &\le C\max\{T(L; \Sigma,t), T(L'; \Sigma, t)\}\max\{\|L\|, \|L'\|\}^{n-1}\|L-L'\|.
    \end{align*}
    In general, we can pick $r\in \Z_{>0}$ divisible enough so that $r\theta, r\theta'$ are represented by $L, L'$ as above. After replacing $t$ by $rt$, we see that both sides are homogeneous of degree $n+1$ in $r$, so we get the same estimate for $\theta, \theta'$. This proves the continuity of $\theta\mapsto \vol(\theta; \Sigma, t)$ on $\Bigc(X)_\Q$ and hence it admits a unique continuous extension to $\Bigc(X)$.

\end{proof}
\begin{remark}
    With a bit more work, it is easy to show that if $s^{-1}\theta'\le \theta\le s\theta'$ for $\theta, \theta'\in \Bigc(X)$ and $s\ge 1$, then for all $t\in \R^d$ with $\min_it_i=0$,
    \[
    s^{-(n+1)}\vol(\theta')\vol(\theta'; \Sigma, st)\le \vol(\theta)\vol(\theta; \Sigma, t)\le s^{n+1}\vol(\theta')\vol(\theta'; \Sigma, s^{-1}t).
    \] 
\end{remark}

For fixed $\theta\in \Bigc(X)$, we record the following concavity property for later use.
\begin{proposition}\label{prop:concave}
For any $\theta\in \Bigc(X)$, the function
    \begin{align*} 
      \vol(\theta; \Sigma, \cdot)\colon  \bbr^{d}&\to \bbr, \quad t\mapsto \vol(\theta;\Sigma,t)
    \end{align*}
    is concave.
\end{proposition}
\begin{proof}
    First assume $L$ is a big line bundle. Let $t_0, t_1\in \R^d$ and denote by $\chi_t$ the geodesic joining $\chi_0:=\chi_{\Sigma,t_0}$ and $\chi_1:=\chi_{\Sigma,t_1}$. For any $m\ge 1$, pick an orthogonal basis $(s_i)$ in $R_m$ for $\chi_0$ and $\chi_1$ and write $s=\sum_i a_i s_i\in R_m$. Then
    \begin{align*}
        \chi_t(s)&=\min_{i,\,a_i\neq 0}\{(1-t)\chi_0(s_i)+t\chi_1(s_i)\}\\
        &=\min_{i,\,a_i\neq 0}\{(1-t)\min_{v\in \Sigma}\{v(s_i)+mt_{0,v}\}+t\min_{v\in\Sigma}\{v(s_i)+mt_{1,v}\}\}\\
        &\leq \min_{i,\,a_i\neq 0}\min_{v\in\Sigma}\{v(s_i)+m((1-t)t_{0,v}+tt_{1,v})\}\\
        &=\min_{i,\,a_i\neq 0}\chi_{\Sigma,(1-t)t_0+tt_1}(s_i)
        \leq \chi_{\Sigma,(1-t)t_0+tt_1}(s).
    \end{align*}
    Combined with~\eqref{eq:affine}, we get
    $$\vol(L;\Sigma,(1-t)t_0+tt_1)\geq \vol(\chi_t)=(1-t)\vol(L;\Sigma,t_0)+t\vol(L;\Sigma,t_1),$$
    where we have used that if $\chi\leq\chi'$ then $\vol(\chi)\leq \vol(\chi')$.
    
    In general, for $\theta\in\Bigc(X)_\Q$, concavity follows from homogeneity, see Lemma~\ref{lem:ratclass}; and the result for $\theta\in \Bigc(X)$ follows from continuity by Proposition~\ref{prop:volcts}.
\end{proof}

\subsection{Extension of Darvas metrics and the Fubini--Study operator to big classes}

Since the $\dd_{1,L}$ pseudometrics on $\cD(\Xlin)$ only depend on the volume function, we may now extend pseudometrics and the Fubini--Study operator on $\cD(\Xlin)$ to big classes.


\begin{definition}\label{def:d1real}
    For $\theta\in\Bigc(X)$, the $\dd_{1, \theta}$ pseudometric and the associated quotient pseudometric on $\cD(\Xlin)$ is given by
    $$\dd_{1,\theta}((\Sigma,t),(\Sigma',t')):=\vol(\theta;\Sigma,t)+\vol(\theta;\Sigma',t')-2\vol(\theta;\Sigma\cup\Sigma',t\wedge t'),$$
    and 
     \[
    \qd_{1,\theta}((\Sigma,t),(\Sigma',t'))\=\inf_{c\in \R}\dd_{1,\theta}((\Sigma,t), (\Sigma',t'+c{\bf 1}))
    \]
    for $(\Sigma,t),(\Sigma',t')\in \cD(\Xlin)$. 
\end{definition}

It follows directly from continuity of volumes that $\dd_{1,\theta}$ are indeed pseudometrics on $\cD(\Xlin)$.
For any $\theta\in\Bigc(X)$, we write $(\cN^\lin(\theta), \de)$ for the Hausdorff quotient of $(\cD(\Xlin), \dd_{1, \theta})$, and similarly for $\cN^\qm(\theta)$.
As a consequence of Theorem~\ref{thm:d1metric} and Corollary~\ref{cor:FSd1}, for $L, L'$ two big line bundles in the same class $\theta$, we get a canonical homeomorphism
\[
(\cN^\lin(\theta), \de)\simeq(\cN^\lin(L), \de)\simeq (\cN^\lin(L'), \de),
\]
which factors through $\rees\colon \cN^\lin(L)\to \cD(\Xlin)$. In particular, the set of Rees valuations $\rees(L;\Sigma, t)$ is independent of the choice of line bundles in the same numerical class.
In view of Corollary~\ref{cor:FSd1}, we may therefore define
\begin{definition}\label{def:FSd1}
    Let $\theta\in \Bigc(X)$ and $(\Sigma, t)\in\cD(\Xlin)$. The set of Rees valuations is
    \[
    \rees(\theta; \Sigma, t)\= \bigcap\{\Sigma'\subseteq \Sigma\mid \dd_{1,\theta} ((\Sigma', t),(\Sigma, t))=0\}.
    \]
    The Fubini--Study operator with respect to $\theta$ is defined as
    $$\FS(\theta;\Sigma,t)(v)\=\inf\{c\in\bbr \mid \dd_{1,\theta}((\Sigma,t),(\Sigma\cup\{v\},t\cup \{c\}))=0\},$$
    for any $v\in \Xlin$.
\end{definition}


\section{Quasimonomial norms and Okounkov bodies}\label{sec:oko}
In this section, we further restrict ourselves to quasimonomial norms and study Okounkov bodies with respect to a quasimonomial valuation. We start by recalling the relevant constructions in~\cite{peng:stab}, and extend his construction to global Okounkov bodies with respect to any numerical class in the big cone.

\subsection{Okounkov bodies with respect to a quasimonomial valuation}
Generalizing a construction in~\cite{lazmus}, Peng~\cite{peng:stab} defines an Okounkov body \emph{with respect to a quasimonomial valuation}. This construction is central to our theory and will be the key ingredient throughout, and so let us briefly recall the construction here. Let $v=v_\alpha\in \QM(Y, E)$ for some log smooth model $(Y, E)$ over $X$ be a quasimonomial valuation on $X$ with $\alpha$ chosen so that all components are $\Q$-linearly independent. Let $L$ be a big line bundle on $X$, and $W_\bullet\subset R(X,L)$ a graded linear series containing an ample series. One can construct as in~\cite[Definition 3.2]{peng:stab} a valuation of maximal rational rank
$$\nu:R(X,L)\to \mathbb{Z}^n$$
and the associated Okounkov body $\Delta(\nu;W_\bullet)$ with respect to this valuation. If $X=Y$ is smooth, and $v$ is divisorial, then the above construction recovers the one in~\cite[Section~4.4]{lazmus}.

We say such a valuation $\nu$ is \emph{adapted to} the quasimonomial valuation, and will sometimes write $\Delta(v; W_\bullet)$ for the Okounkov body $\Delta(\nu; W_\bullet)$, when $\nu$ is adapted to $v\in \Xqm$. We also say that $\Delta(v; W_\bullet)$ is an Okounkov body \emph{with respect to} $v$. Such an Okounkov body evidently depends on the choice of $\nu$, but for our purposes, we will always fix a choice of $\nu$ as long as it is adapted to $v$.

Denote by $V^\lambda_\bullet(W_\bullet, \{v\}, \{t_v\})$ the graded linear series induced by the norm $\chi\=v+t_v$ on $W_\bullet$ for some $t_v\in \R$, we record the following version of a result in~\cite{peng:stab}.
\begin{theorem}[{\cite[Theorem~1.3, Proposition~3.13]{peng:stab}}]\label{thm:qmokbody}
    Let $W_\bullet$ be a graded series of $R(X,L)$ containing an ample series and $v\in \Xqm$ be a quasimonomial valuation of rational rank $r$. For $\lambda\in [t_v, T(W_\bullet; \{v\}, \{t_v\}))$, we have
    \begin{equation}
        \Delta(\nu;V_\bullet^\lambda(W_\bullet;\{v\},\{t_v\}))=\Delta(\nu; W_\bullet)\cap \{x\in \R^n\mid \langle \alpha, x\rangle+t_v \ge \lambda\}
    \end{equation}
    and thus
    \begin{equation}
        \vol(V^\lambda_\bullet(W_\bullet;\{v\},\{t_v\}))=n!\vol_{\bbr^n}(\Delta(\nu; W_\bullet)\cap \{x\in \R^n\mid \langle \alpha, x\rangle+t_v \ge \lambda\}),
    \end{equation}
    where by abuse of notation, we write $\langle \alpha, x\rangle$ for the usual inner product $\langle \alpha, \pi(x)\rangle$ on $\R^r$, and $\pi:\bbr^n\to\bbr^{r}$ is the projection onto the first $r$-coordinates.
    
\end{theorem}

Partly motivated by~\cite[Lemma~3.28]{yli}, we can use Okounkov bodies with respect to a quasimonomial valuation to compute the value of the Fubini--Study operator at the given valuation.
\begin{proposition}\label{prop:formula_concave}
    Let $\chi$ be a norm on $R\=R(X, L)$ for $L$ big, and $v\in X^\qm$ a valuation of rational rank $r$. Then for an Okounkov body $\Delta(v; L)$ with respect to $v$, we have
    $$\FS(\chi)(v)=\sup_{p\in\Delta(v;L)^\circ}\left(G[\chi](p)-\langle \alpha,p\rangle\right).$$
\end{proposition}
\begin{proof}
    Let $\nu$ be the valuation adapted to $v$ giving rise to the Okounkov body $\Delta(v; L)$. Then by construction we have 
    $\langle \alpha,\nu(s)\rangle=v(s)$ for any $s\in R_m\setminus\{0\}$. Now, for the finite-dimensional concave transform in degree $m$, and $p\in \Gamma_m(L)$, we have
    $$G_m[\chi](p)=\sup\{\chi(s)\mid s\in R_m,\,\nu(s)=p\},$$
    and thus
    \[
    \frac 1m\max_{p\in \Gamma_m(L)} (G_m[\chi](p)-\langle \alpha,p\rangle) = \frac 1m\sup_{s\in R_m\setminus\{0\}}(\chi(s)-v(s)).
    \]
    Taking supremum over $m$, we get
    \[
    \FS(\chi)(v)=\sup_m\max_{p\in \Gamma_m(L)} \left(\frac 1m (G_m[\chi](p)-\langle \alpha,p\rangle)\right)
    \]
    Note also that for each $p\in \Delta(L)^\circ$, 
    \[
    G[\chi](p)=\sup_{p_i\in \Gamma_{m_i}(L), p_i/m_i\to p}\frac{1}{m_i}G_{m_i}[\chi](p_i),
    \]
    after taking supremum over $p\in \Delta(L)^\circ$, we have
    \[
    \FS(\chi)(v)=\sup_{p\in\Delta(L)^\circ}\left(G[\chi](p)-\langle \alpha,p\rangle\right).
    \]
\end{proof}

Restricting to valuative norms, we note the following elementary but important induction procedure. Suppose $\chi=\chi_{\Sigma,t}$ is a valuative norm supported on $\Sigma=\{v_i\}_{i=1}^d$. For each $i=1, \dots, d$, let us write 
\[
V_{i,\bullet}^\lambda(L;\Sigma,t)\= V_\bullet^\lambda(L;\Sigma\setminus\{v_i\},t\setminus\{t_i\}),
\]
where if $d=1$, $V_{1, \bullet}^\lambda\=R(X,L)$ for any $\lambda$.
Then for $\lambda<T(L;\Sigma,t)$ and each $i=1,\dots,d$, we have
\begin{equation}\label{eq:indseries}
    V^\lambda_\bullet(L;\Sigma,t)
    = V_\bullet^\lambda(V_{i,\bullet}^\lambda(L;\Sigma,t);\{v_i\},t_i)
\end{equation}
As a direct corollary of Theorem~\ref{thm:qmokbody}, the inductive process gives the following characterization of concave transforms on an Okounkov body with respect to $v\in\Xqm$.
\begin{corollary}\label{cor:concavetr}
    Suppose $|\Sigma|\geq 2$. Let $\chi\=\chi_{\Sigma, t}, \chi'\=\chi_{\Sigma\setminus\{v\}, t\setminus\{t_v\}}$ be valuative norms with $v\in \Xqm$, and $\Delta(v; L)$ be an Okounkov body with respect to $v$. Then
    \[
    G[\chi](p)=\min\{G[\chi'](p), \langle \alpha, p\rangle+t_v\}
    \]
    for any $p\in \Delta(v; L)$.
\end{corollary}
\begin{proof}
    This follows from unwinding definitions. Indeed, by Theorem~\ref{thm:qmokbody} we have for any $\lambda<T(L; \Sigma, t)$,
    \[
    \Delta(v;V_\bullet^\lambda(\chi))=\Delta(v; V_\bullet^\lambda(\chi'))\cap \{\langle \alpha, x\rangle+t_v \ge \lambda\}
    \]
    Thus for any $p\in \Delta(v; L)$,
    \begin{align*}
    G[\chi](p)&=\sup\{\lambda\mid p\in \Delta(v; V_\bullet^\lambda(\chi'))\cap \{\langle \alpha, x\rangle+t_v \ge \lambda\}\}\\
    &=\min\{G[\chi'](p), \sup\{\lambda\mid p\in \{\langle\alpha, x\rangle+t_v\ge \lambda\}\}\}\\
    &=\min\{G[\chi'](p), \langle \alpha, p\rangle+t_v\}.
    \end{align*}
\end{proof}

\subsection{Global Okounkov bodies}\label{sec:globalok}
We now extend the construction of an Okounkov body with respect to a quasimonomial valuation to a global version for numerical classes in the big cone, following~\cite[\S4]{lazmus}. 
As in~\cite[\S4.2]{lazmus}, we pick a $\Z$-basis $L_1, \dots, L_r$ of $N^1(X)$ such that the pseudoeffective cone is contained in the positive orthant of $\R^r$. For each $\vec{\gamma}=(\gamma_1, \dots, \gamma_r)\in \N^r$, we set $\vec{\gamma}L\=\sum_{i=1}^r\gamma_iL_i$ and $|\vec{\gamma}|=\sum_i\gamma_i$.
\begin{proposition}\label{prop:bigcvx}
    Denote by $S\=\{\vec\gamma L,\,\gamma_i>0\,\forall i,\,|\vec\gamma|=1\}\cap \Bigc(X)$. Let $\Sigma=\{v_i\}_{i=1}^d\subset X^\qm$ be a finite set and $\nu$ a rank $n$ valuation adapted to some $v\in \Xqm$. There exists a closed convex set
    $$\Delta(\nu;\Sigma)_S\subset S\times \bbr\times \bbr^{|\Sigma|}\times\bbr^n,$$
    whose fibers $\Delta^\lambda(\nu;\theta;\Sigma,t)\subset \bbr^n$ over $(\theta,\lambda,t)$ are either empty or convex bodies, all contained in a fixed compact subset of $\bbr^n$.
    
    If $\theta=c_1(L)$ with $L$ a big $\bbq$-line bundle, then  for any $t$ and any $\lambda< T(L;\Sigma,t)$,
    \begin{equation}\label{eq:identification}
        \Delta^\lambda(\nu;\theta;\Sigma,t)=\Delta(\nu;V^\lambda_\bullet(L;\Sigma,t)).
    \end{equation}
    In particular, if $\nu=\nu_i$ is adapted to some $v_i\in \Sigma$, then for     
    $\lambda<T(L;\Sigma,t)$, we have
    $$\Delta^\lambda(\nu_i;\theta;\Sigma,t)=\Delta(\nu_i;V^\lambda_\bullet(L;\Sigma\setminus\{v_i\},t\setminus\{t_i\}))\cap \{\langle \alpha_i,p\rangle+t_i\geq \lambda\}.$$
\end{proposition}
\begin{proof}
    We consider the $(\mathbb{Z}^r\times \mathbb{Z}\times \mathbb{Z}^{|\Sigma|})$-graded  linear series
    $$V_{\vec\gamma, \lambda,t}:=\{s\in H^0(X,\vec \gamma L),\,\min_i(v_i(s) + t_i) \geq \lambda\}$$
    for $\vec\gamma\in\mathbb{Z}^r$, $\lambda\in\mathbb{Z}$, $t\in \mathbb{Z}^{|\Sigma|}$.
    Define $\Delta(\nu;\Sigma)_S$ to be the intersection with $(S\times\bbr\times\bbr^{|\Sigma|}\times\bbr^n)$ of the closure
    \begin{equation}
        \overline{\left(|\vec\gamma|\mi(\vec\gamma L,\lambda,t,\nu(s)),\,\vec\gamma\in\mathbb{Z}^r\setminus\{0\},\,\lambda\in\mathbb{Z},\,t\in \mathbb{Z}^{|\Sigma|},\,s\in V_{\vec\gamma, \lambda,t}\setminus\{0\}\right)}
    \end{equation}
    That $\Delta(\nu; \Sigma)_S$ is bounded follows from the standard uniform linear bound on each coordinate $|\nu(s)_i|\leq C|\vec\gamma|$ for $s\in H^0(X,\vec\gamma L)$ restricted to $|\vec\gamma|=1$, see~\cite[Theorem~3.9]{peng:stab}. It is a convex set as follows from the grading of $V$, since intersecting with the convex slice $S$ preserves convexity. This proves the first assertion.

    We now prove \eqref{eq:identification}. We begin with $\gamma_i,\lambda,t_i\in\mathbb{Z}$ such that $L:=\sum_i \gamma_i L_i$ is big. Consider the multigraded linear series $W_{\vec \bullet}^\lambda$ with
    \[
        W_{\vec\gamma}^\lambda\=V_{\vec\gamma,|\vec\gamma|\lambda,|\vec\gamma|t}.
    \]
    It is easy to check that $\mathrm{supp}(W_{\vec\bullet}^\lambda)$ has non-empty interior if $\lambda<T(\theta;\Sigma, t)$, and that $W^\lambda_{\vec\bullet}$ contains an ample series for some $\vec\gamma\in \mathrm{int}(\mathrm{supp}(W^\lambda_{\vec\bullet}))$, in the sense of~\cite[Definition~4.17, Theorem~4.21]{lazmus}. Using \cite[Corollary A.3]{lazmus} applied to the semigroup $\Gamma:=\{(\nu(s),\vec\gamma,\lambda,t),\,s\in V_{\vec\gamma,\lambda,t}-\{0\}\}$, the projection $p$ omitting the $p$-coordinate in $\bbr^n$, and the linear subspace $L_{\lambda,t}:=\{\vec\gamma,(\sum_i \gamma_i)\lambda,(\sum_i \gamma_i)t\}$, one has
    $$\overline{\mathrm{cone}}(\Gamma)\cap p\mi(L_{\lambda,t})=\overline{\mathrm{cone}}(\Gamma\cap p\mi(L_{\lambda,t}))=\overline{\mathrm{cone}}(\{(\nu(s),\vec\gamma,|\vec\gamma|\lambda,|\vec\gamma|t),s\in W^\lambda_{\vec\gamma}-\{0\}\}).$$
    Applying \cite[Proposition 4.9]{lazmus} to the fibre over $\theta=c_1(L)/|\vec\gamma|$ then gives
    \[
    \Delta^\lambda(\nu; \theta; \Sigma, t)=|\vec\gamma|\mi\Delta(\nu; V_\bullet^{|\vec\gamma|\lambda}(L; \Sigma, |\vec\gamma|t)).
    \]
    For rational $\vec\gamma,\lambda,t$, one extends by homogeneity using
    $$V_m^{q\lambda}(qL;\Sigma,qt)=V^\lambda_{qm}(L;\Sigma,t).$$
    Since $\Delta(\nu;\Sigma)_S$ is convex and its fibres for $\lambda<T(L;\Sigma,t)$ are nonempty, those fibres are continuous in $(t,\lambda)$ in the Hausdorff distance. By continuity of both sides, this extends \eqref{eq:identification} to real $\lambda,t$. 
    

\end{proof}

\begin{remark}
    Fixing $t$ and $\theta$, the family $\lambda\mapsto \Delta^\lambda(\nu;\theta;\Sigma,t)$ is an Okounkov test curve in the sense of \cite{xia:oko}. A similar construction, for a fixed $\Sigma=\{\nu\}$ containing a single valuation, and without the $(\lambda,t)$-coordinates, appears in \cite[Section 3.3.4]{peng:stab}.
\end{remark}
Going forward, we will write
\begin{equation}
\Delta(\nu;\Sigma):=\bbr_{>0} \Delta(\nu;\Sigma)_S=\{(x\theta,x\lambda,xt,xp)\mid x>0,\,(\theta,\lambda,t,p)\in \Delta(\nu;\Sigma)_S\}
\end{equation}
for the cone over $\Delta(\nu;\Sigma)_S$, and use the same notation for fibers as in the previous proposition. 

\subsection{Filtered Okounkov bodies}
We next generalize the filtered Okounkov bodies in~\cite{bouchen} to our setting.
\begin{definition}\label{def:filtered}
    The \textit{filtered Okounkov body} is
    $$\hat \Delta(v;\theta;\Sigma,t):=\{(\lambda,p)\mid(\theta;\lambda,t,p)\in \Delta(v;\Sigma),\,\lambda\geq \min_i t_i\}\subset \bbr\times\bbr^n.$$
\end{definition}
When $\theta=c_1(L)$,
this coincides by construction, with the one defined in~\cite[Definition~1.8]{bouchen}, that is,
\begin{equation}\label{eq:filtered}
    \hat\Delta(v;\theta;\Sigma,t)=\{(x,p)\mid p\in \Delta(\nu;L),\,\min_i t_i\leq x \leq G[\theta;\Sigma,t](p)\},
\end{equation}
    where $G[\theta; \Sigma, t]=G[\chi_{\Sigma, t}]$ is the concave transform.
Thus $\Delta(v;\Sigma)$ encodes the information about concave transforms of all quasimonomial data supported on $\Sigma$, with respect to all big classes.
In particular, it is easy to see from~\eqref{eq:volint} and~\cite[Definition 1.9, Corollary 1.13]{bouchen} that
\begin{equation}\label{eq:volcomp}
    \frac{n!}{\vol(\theta)}\vol_{\bbr^{n+1}}(\hat \Delta(v;\theta;\Sigma,t))=\vol(\theta;\Sigma,t)-\min_i t_i.
\end{equation}
For truncations by any $c\leq \min_i t_i$ in the definition of filtered Okounkov bodies, which we will simply call "filtered Okounkov bodies truncated at $c$", and denote in the same way, we obtain similarly
\begin{equation}\label{eq:volcomp2}
    \frac{n!}{\vol(\theta)}\vol_{\bbr^{n+1}}(\hat \Delta(v;\theta;\Sigma,t))=\vol(\theta;\Sigma,t)-c.
\end{equation}

The following is a generalization of Proposition~\ref{prop:formula_concave}.
\begin{proposition}\label{prop:fs_formula_real}
    Let $\theta\in \Bigc(X)$, $v=v_\alpha\in X^\qm$ with adapted flag valuation $\nu$, $\Sigma\subset X^\qm$ finite, $t\in \bbr^{|\Sigma|}$. Then
    $$\FS(\theta;\Sigma,t)(v)=\sup_{(\lambda,p)\in\hat\Delta(v;\theta;\Sigma,t)}\left(\lambda - \langle \alpha,p\rangle\right).$$
\end{proposition}
\begin{proof}
    Without loss of generality we may assume $\min_i t_i=0$. Let $\nu$ be a valuation adapted to $v=v_\alpha$. 
    Suppose $c$ is a candidate in the definition of $\FS$, and note that in particular
    \begin{equation}
        \de(\theta;(\Sigma,t),(\Sigma\cup\{v\},t\cup \{c\}))=\vol(\theta;\Sigma,t)-\vol(\theta,\Sigma\cup\{v\},t\cup \{c\})=0.
    \end{equation}
    By \eqref{eq:volcomp2} (for filtered Okounkov bodies truncated at $\min(0,c)$) this gives
    \begin{equation}
        \vol_{\bbr^{n+1}}(\hat\Delta(v;\theta;\Sigma,t))-\vol_{\bbr^{n+1}}(\hat\Delta(v;\theta;\Sigma,t)\cap\{\langle \alpha,p\rangle+c\geq \lambda\})=0,
    \end{equation}
    where we have used that the valuation used to compute $\hat\Delta$ is adapted to $v$ (this holds for real classes as well by approximating by rational $\theta_k$ with $\theta_k-\theta$ ample). Since $\hat\Delta(v;\theta;\Sigma,t)$ is a convex body and $\{\langle \alpha,\pi(p)\rangle+c\geq \lambda\})$ is a closed half-space, the above volume equality is equivalent to
    $$\hat\Delta(v;\theta;\Sigma,t)\cap\{\langle \alpha,p\rangle+c\geq \lambda\})=\hat\Delta(v;\theta;\Sigma,t),$$
    that is, $\lambda-\langle \alpha,p\rangle\leq c$ for all $(\lambda,p)\in \hat\Delta(v;\theta;\Sigma,t)$. This concludes the proof.
\end{proof}

\begin{proposition}\label{prop:d1uniqueness}
    For $t,t'\in \bbr^{|\Sigma|}$, and $\alpha\in |0,1]$, if $\theta$ is a big real class, and
    $$\alpha\mapsto \vol(\theta;\Sigma,(1-\alpha)t+\alpha t')$$
    is affine, then
    \[
    \qd_{1,\theta}((\Sigma,t), (\Sigma,t'+c{\bf 1}))=\inf_{c\in \R}\dd_{1,\theta}((\Sigma,t), (\Sigma,t'+c{\bf 1}))=0,
    \]
    and the infimum is achieved by $c=\vol(\theta; \Sigma, t')-\vol(\theta; \Sigma, t)$.
\end{proposition}
\begin{proof}
    Wlog setting $t\to t-\vol(L;\Sigma,t),t'\to t'-\vol(L;\Sigma,t)$ we can assume that $\vol(L;\Sigma,t)=\vol(L;\Sigma,t')=0$, hence that
    \begin{equation}\label{eq:normalisation}
    \vol(L;\Sigma,(1-\alpha)t+\alpha t')=0.
    \end{equation}
    Pick now some $i$ and let $\nu_i$ be adapted to the qm valuation $v_i=v_{\alpha_i}$. Pick $\alpha\in[0,1]$, and we will consider filtered bodies $\hat\Delta$ in the sense of Definition \ref{def:filtered}. By Proposition \ref{prop:bigcvx}, using convexity of the global body $\Delta(\nu;\Sigma)$, one obtains
    \begin{equation}\label{eq:incl}
        \hat\Delta(\nu_i;\theta;\Sigma,(1-\alpha)t+\alpha t')\supseteq(1-\alpha)\hat\Delta(\nu_i;\theta;\Sigma,t)+\alpha \hat\Delta(\nu_i;\theta;\Sigma,t').
    \end{equation}
    Passing to volumes, and using the Brunn--Minkowski inequality,
    \begin{equation}\label{eq:vol}
    \begin{aligned}
       & \vol_{\bbr^{n+1}}\left(\hat\Delta(\nu_i;\theta;\Sigma,(1-\alpha)t+\alpha t')\right)\\
        &\geq \vol_{\bbr^{n+1}}\left((1-\alpha)\hat\Delta(\nu_i;\theta;\Sigma,t)+\alpha \hat\Delta(\nu_i;\theta;\Sigma,t')\right)\\
        &\geq \left((1-\alpha)\vol_{\bbr^{n+1}}(\hat\Delta(\nu_i;\theta;\Sigma,t))^{\frac1{n+1}}+\alpha\vol_{\bbr^{n+1}}(\hat\Delta(\nu_i;\theta;\Sigma,t'))^{\frac1{n+1}}\right)^{n+1}.
    \end{aligned}
    \end{equation}
    Thus \eqref{eq:normalisation} and \eqref{eq:volcomp2} (applied to filtered bodies truncated at $\min(0,\min_i t_i,\min_i t'_i)$) implies equality in \eqref{eq:vol}. By the equality case of Brunn--Minkowski, we then have that $\hat\Delta(\nu_i;\theta;\Sigma,t)$ and $\hat\Delta(\nu_i;\theta;\Sigma,t')$ are homothetic. Since they have the same volume by assumption, they must in fact be translates. Finally, normalisation (for the $\lambda$-component of the translation) and the fact that they have the same projection to the $p$-coordinates implies the two bodies to be equal. Using Proposition \ref{prop:fs_formula_real} we then get
    \begin{equation}
        \FS(\theta;\Sigma,t)(v_i)=\FS(\theta;\Sigma,t')(v_i)
    \end{equation}
    for all $i$. By the expression for $\FS$ in Definition \ref{def:FSd1} (applied to each coordinate successively), we then have
    $$\dd_{1,\theta}((\Sigma,t),(\Sigma,\FS(\theta;\Sigma,t)))=\dd_{1,\theta}((\Sigma,t'),(\Sigma,\FS(\theta;\Sigma,t)))=0,$$
    hence
    $$\dd_{1,\theta}((\Sigma,t),(\Sigma,t'))=0$$
    by triangle inequality. This concludes the proof.    
\end{proof}

\section{Quasimonomial Monge--Ampère equations}\label{sec:MA}

In this section, we define and study the Monge--Amp\`ere operator on quasimonomial norms with respect to a big class, and prove the Calabi--Yau theorem for quasimonomial measures.

\subsection{The Monge--Ampère operator with respect to a big line bundle}
Let $L$ be a big line bundle and $\theta\=c_1(L)$. We first define the Monge--Amp\`ere operator 
\[
\MA_\theta\colon \cD(\Xqm)\to \cM^\qm,
\]
where $\cM^{\qm}$ is the set of \emph{quasimonomial measures}, i.e. probability measures on $\Xna$ supported on finite subsets of $\Xqm$. 
The key result for making our definition is the following differentiability result for the volume function.
\begin{theorem}\label{thm:diff}
    Let $L$ be a big line bundle, and $\Sigma\=\{v_i\}_{i=1}^d$ be a fixed set of quasimonomial valuations. Then $\vol(L; \Sigma, \cdot)\in C^1(\R^d)$. Moreover, 
    \begin{enumerate}[(i), font=\normalfont]
        \item $\nabla \vol(L; \Sigma, t)\in \sigma\=\{\xi=(\xi_1, \dots, \xi_d)\in \R^d_{\ge0}\mid \sum_{i=1}^d\xi_i=1\}$;
        \item   $\partial_i\vol(L; \Sigma, t)>0$ iff $v_i\in \rees(\chi_{\Sigma, t})$, and $\partial_i\vol(L; \Sigma, t)=0$ otherwise.
    \end{enumerate}
\end{theorem}
We need the following elementary lemma.
\begin{lemma}
    Let $\Delta\subset\R^n$ be a nonempty compact convex set and $G$ a concave function on $\Delta$. For any $c\in \R$ with $\{G>c\}\subset \Delta$ nonempty, we have
    \[
    \rho(\{G>c\})=\rho(\{G\ge c\}),
    \]
    with $\rho$ the Lebesgue measure on $\Delta$.
\end{lemma}
\begin{proof}
    This follows from concavity of $G$. Indeed, for any $x\in\{G>c\}$ and $y\in \{G\ge c\}$, we have $(1-t)x+ty\in\{G>c\}$ for $t\in [0,1)$. Thus for fixed $x\in \{G> c\}$, we have
    \[
    t^n\rho(\{G\ge c\})=\rho (t\{G\ge c\})\le \rho(\{G>c\})\le \rho(\{G\ge c\}),
    \]
    for any $t\in [0,1)$. We conclude by letting $t\to 1$.
\end{proof}
\begin{proof}[Proof of Theorem~\ref{thm:diff}]
    To simplify the notation, we write $\chi(t)\=\chi_{\Sigma, t}$, $\chi_i(t)\=\chi_{\Sigma\setminus\{v_i\}, t\setminus\{t_i\}}$, and $\Delta_i\=\Delta(\nu_i; L)$ for the Okounkov body with respect to $v_i=v_{\alpha_i}$ for $i=1, \dots, d$. Let $e_i\in \R^d$ be the $i$-th coordinate vector. We also write $$\partial_i^\pm\vol(L; \Sigma, t)\=\frac{d}{d\varepsilon}\Big|_{\varepsilon=0^\pm}\vol(L; \Sigma, t+\varepsilon e_i),$$
    provided that it exists.

    We first show all partial derivatives exist. Assume first that $\Sigma=\rees(\chi(t))$ and fix some $i\in\{1, \dots, d\}$. Then by Corollary~\ref{cor:concavetr}, we have
    \[
    G[\chi(t)](p)=\min\{G[\chi_i(t)](p), \langle\alpha_i, p\rangle+t_i\},
    \]
    for any $t\in \R^d, p\in \Delta_i$. From~\eqref{eq:volint} and dominated convergence, we deduce
    \[
    \partial_i^+\vol(L; \Sigma, t)
    =\frac{1}{\vol(\Delta_i)}\vol(\{p\in \Delta_i\mid G[\chi_i(t)](p)>\langle\alpha_i, p\rangle+t_i\}),
    \]
    and
    \[
    \partial_i^-\vol(L; \Sigma, t)
    =\frac{1}{\vol(\Delta_i)}\vol(\{p\in \Delta_i\mid G[\chi_i(t)](p)\ge\langle\alpha_i, p\rangle+t_i\}).
    \]    
    We need to show that the left and right derivatives agree. Since $v_i\in \rees(\chi(t))$, there exists $s\in R_m$ with $v_i(s)+mt_i<\chi_i(t)(s)$, giving 
    \[
    t_i< \frac 1m (\chi_i(t)(s)-v_i(s))\le \FS(\chi_i)(v_i)= \sup_{p\in \Delta_i^\circ}(G[\chi_i(t)](p)-\langle \alpha_i, p\rangle),
    \]
    where the last equality follows from Proposition~\ref{prop:formula_concave}. In particular, $\{G[\chi_i(t)]>\langle\alpha_i, p\rangle+t_i\}\ne \emptyset$, and the previous lemma applies. In this case we have
    \[
    \partial_i\vol(L; \Sigma, t)=\frac{1}{\vol(\Delta_i)}\vol(\{p\in \Delta_i\mid G[\chi_i(t)](p)>\langle\alpha_i, p\rangle+t_i\})>0.
    \]
    Moreover, note that for any $c\in \R$, $\vol(\chi(t)+c{\bf 1})=\vol(\chi(t))+c$. By Proposition~\ref{prop:concave} we have
    \begin{equation}\label{eq:probmeas}
    \pm 1=\vol(L; \Sigma, t\pm {\bf 1})-\vol(L; \Sigma, t)\le \pm\sum_{j=1}^d \partial_j\vol(L; \Sigma, t),
    \end{equation}
    and so $\sum_{j=1}^d \partial_j\vol(L; \Sigma, t)=1$.

    In general, wlog we may assume $\Sigma\setminus\{v_i\}=\rees(\chi(t))=\rees(\chi_i(t))$. In particular, $\vol(\chi(t))=\vol(\chi_i(t))$ and
    $\partial_j\vol(\chi_i(t))=\partial_j\vol(L; \Sigma,t)$ for $j\ne i$. We claim that $\partial_i\vol(L; \Sigma, t)=0$. Indeed,~\eqref{eq:probmeas} applied to a supergradient $p=(p_j)_{j=1}^d$ now gives $p_i\le 0$; on the other hand, $\vol(L; \Sigma, t)$ is increasing in each variable, and so we must have $p_i=0$. This proves~(i) and~(ii).

    It remains to show all partial derivatives $\partial_i\vol(L; \Sigma, t)$ are continuous in $t$. By a slight abuse of notation, let us write $\nabla\vol(L; \Sigma, t)\=(\partial_i\vol(L; \Sigma, t))_{i=1}^d\in \R^d$.
    Suppose $t_k\to t\in \R^d$. Then by what precedes, we have that all $\nabla\vol(L; \Sigma, t_k)$ lie in the compact set $\sigma$. After passing to a subsequence, we may assume $\nabla\vol(L; \Sigma, t_k)\to \xi$ for some $\xi\in\sigma$. By concavity of the volume function again, for any $t'\in \R^d$ we have
    \[
    \vol(L; \Sigma, t')\le \vol(L; \Sigma, t_k)+\langle \nabla\vol(L; \Sigma, t_k), t'-t_k\rangle\to \vol(L; \Sigma, t)+\langle \xi, t'-t\rangle,
    \]
    as $k\to\infty$. This shows $\xi$ is a supergradient at $t$, and so $\xi=\nabla\vol(L; \Sigma, t)$. 
    
    

\end{proof}

As a direct consequence, we get the following explicit expression for the partial derivatives of the volume function.
\begin{corollary}\label{cor:partialMA}
    In the notation of Theorem~\ref{thm:diff} and fix a rank $n$ valuation $\nu_i$ adapted to $v_i$, we have for any $i=1,\dots, d$,
    \[
    \partial_i\vol(L; \Sigma, t)
    =\frac{n!}{\norm[\alpha_i]\vol(L)}\int_\bbr \vol_{\bbr^{n-1}}(\Delta(\nu_i; V^\lambda_{i, \bullet}(L;\Sigma,t))\cap \{x\in \R^n\mid \langle\alpha_i, x\rangle +t_i=\lambda\})\,d\lambda.
    \]
\end{corollary}
\begin{proof}
    By the proof of Theorem~\ref{thm:diff} and using the notation in the proof, we have for each $i$, 
    \[
    \partial_i\vol(L; \Sigma, t)=\frac{n!}{\vol(L)}\vol(\{p\in \Delta(\nu_i;L)\mid\,G[\chi_i(t)](p)\geq \langle \alpha_i,p\rangle+t_i\}).
    \]
    On the other hand, since $G$ is the concave transform, for any $\lambda<T(L; \Sigma, t)$, we have
    \[
    \{G[\chi_i(t)]\geq \lambda\}=\Delta(\nu_i; V^\lambda_{i,\bullet}(\Sigma, t)),
    \]
    for $V_{i, \bullet}^\lambda$ as defined in~\eqref{eq:indseries}.
    The result follows by applying the coarea formula to the above set and the function $p\mapsto \langle \alpha_i,p\rangle + t_i$.
\end{proof}

\begin{remark}
    Analogous to~\cite[Theorem~4.26]{lazmus}, the integrand in the previous corollary may be thought of as a generalization of the restricted volume function to a quasimonomial valuation.
\end{remark}

Recall from Lemma~\ref{lem:ratclass} that $\vol(L; \Sigma, t)$ only depends on the numerical class, differentiability of the volume function now allows us to define the Monge--Ampère operator with respect to an integral class.
\begin{definition}
    Let $L$ be a big line bundle and $\theta\=c_1(L)$. The \emph{Monge--Amp\`ere operator} is the map
    \[
    \MA_\theta\colon \cD(\Xqm)\to \cM^\qm, (\Sigma=\{v_i\}_{i=1}^d, t)\mapsto \MA_\theta(\Sigma, t) \= \sum_{i=1}^d\MA_\theta(\Sigma, t)_i\delta_{v_i},
    \]    
    where $\MA_\theta(\Sigma,t)_i\=\partial_i\vol(L; \Sigma, t)$ is the $i$-th component of $\nabla\vol(L; \Sigma, t)$.
    By abuse of notation, we will also denote by $\MA_\theta(\Sigma, t)\=\nabla\vol(\theta; \Sigma, t)\in \sigma_\Sigma$, where $\sigma_\Sigma$ is the standard simplex in $\R^d$.
\end{definition}
By Theorem~\ref{thm:diff}, the Monge--Amp\`ere measure $\MA_\theta(\Sigma, t)$ is a translation invariant quasimonomial measure supported on $\rees(\chi_{\Sigma, t})$.
The next proposition shows that the above map factors through $\cN^\qm(L)$ to define the \emph{Monge--Amp\`ere operator} on quasimonomial norms
\[
\MA_\theta\colon \cN^\qm(L)\to \cM^\qm, \chi_{\Sigma, t}\mapsto \MA_\theta(\Sigma, t).
\]
\begin{proposition}\label{prop:propertiesMA}
    Suppose $(\Sigma=\{v_i\}_{i=1}^d,t), (\Sigma', t')\in\cD(\Xqm)$ satisfies $\chi_{\Sigma, t}=\chi_{\Sigma', t'}\in \cN^\qm(L)$ for a big line bundle $L$. Then
    $$\MA_L(\Sigma,t)=\MA_L(\Sigma',t').$$
\end{proposition}
\begin{proof}
    First assume $\min_{i=1, \dots, d} t_i=0$ using translation invariance. Note that if $\Sigma\subset\Sigma'$ and $t_v>\max_{v_j\in\Sigma} T_L(v_j)$ for $v\in \Sigma'\setminus\Sigma$, then in the notation of~\eqref{eq:embedding}, $v\notin\rees(\chi_{\Sigma', \iota(t)})$, and by Theorem~\ref{thm:diff} we have that the coefficients of $\MA_L(\Sigma',\iota(t))$ corresponding $v\in \Sigma'\setminus\Sigma$ must be zero. Thus after replacing $\Sigma, \Sigma'$ by $\Sigma\cup \Sigma'$, and $t,t'$ by $\iota(t), \iota(t')$, we may assume $\Sigma=\Sigma'$.
    We first show that if $\chi_{\Sigma,t}=\chi_{\Sigma,t'}$ with $t_i\geq t'_i$ for all $i$, then $\MA_L(\Sigma,t)=\MA_L(\Sigma,t')$. Indeed, consider for each coordinate vector $e_i\in\bbr^{d}$ the functions
    $$f(\varepsilon):=\vol(L;\Sigma,t+\varepsilon e_i),\quad g(\varepsilon):=\vol(L;\Sigma,t'+\varepsilon e_i).$$
    By assumption we have $f(\varepsilon)\geq g(\varepsilon)$, $f(0)=g(0)$. Since both $f, g\in C^1(\R^d)$ by Theorem~\ref{thm:diff}, we get $f'(0)=g'(0)$ by taking $\varepsilon\to 0+$ and $\varepsilon\to 0-$.

    In general, by Lemma~\ref{lem:lattice} we may then find $s\le t\wedge t'$ coordinatewise (which can be taken to be the minimal $t_0$ in Proposition \ref{prop:reesnorm}(i)) with 
    \[
    \chi_{\Sigma, s}=\chi_{\Sigma, t}=\chi_{\Sigma, t'},
    \]
    and we conclude by the previous case.
\end{proof}

We record here the simpler expression for the Monge--Ampère measure in the divisorial case.
\begin{proposition}
    Let $L$ be a big line bundle on $X$, and $(\Sigma, t)$ be a divisorial datum, with each $v_j=c_j\ord_{E_j}\in \Sigma$. Let $e_j$ be the $j$-th unit vector in $\R^d$. Let $\pi: Y\to X$ be a smooth model with all $E_i$ being prime divisors on $Y$.
    Then
    \[
    \frac{d}{d\varepsilon}\Bigg|_{\varepsilon=0}\vol(L;\Sigma,t+\varepsilon e_j)=\frac{n}{\vol(L)}c_j\mi\int_{t_j}^\infty\vol_{Y|E_i}(\pi^*L-\sum_{i=1}^dc_i^{-1}(\lambda-t_i)_+E_i).
    \]
    Hence
    \[
    \MA_L(\Sigma,t):=\frac{n}{\vol(L)}\sum_{j}\left(c_j\mi\int_{t_j}^\infty\vol_{Y|E_j}(\pi^*L-\sum_{i=1}^dc_i^{-1}(\lambda-t_i)_+E_i)\,d\lambda\right)\delta_{v_j}.
    \]
\end{proposition}
\begin{proof}
    Using Corollary \ref{cor:partialMA}, it will suffice to show
    \begin{align*}
    &(n-1)!\vol(\Delta^0(\nu_i; V^\lambda_\bullet(L;\Sigma\setminus\{v_i\},t\setminus\{t_i\}))\cap \{(x_1, \dots, x_d)\in \R^d\mid c_ix_1+t_i=\lambda\})\\
    &=\vol_{Y|E_i}(\pi^*L-\sum_{i=1}^dc_i^{-1}(\lambda-t_i)_+E_i),
    \end{align*}
    which follows from~\cite[Proposition~4.7]{bj:ytd} and~\cite[Corollary C]{bfj:teissier}.
\end{proof}

\subsection{Extension of the Monge--Ampère operator to real classes}
We now extend the Monge--Amp\`ere operator to arbitrary $\theta\in \Bigc(X)$. By what precedes, it amounts to showing that if $\theta\in \Bigc(X)$ then
$$t\mapsto \vol(\theta;\Sigma,t)$$
is also differentiable. For $\theta\in\Bigc(X)_\Q$, this readily follows from the homogeneity of the volume, see Lemma~\ref{lem:ratclass}. In general, we use the global Okounkov body construction described in~\S\ref{sec:globalok}.
\begin{definition}
    For $\theta\in \Bigc(X)$, the \emph{Monge--Amp\`ere operator} with respect to $\theta$ is defined as
    \[
    \MA_\theta\colon \cD(\Xqm) \to \cM^\qm, (\Sigma=\{v_i\}_{i=1}^d, t)\mapsto \MA_\theta(\Sigma, t)\=\sum_{i=1}^d\MA_\theta(\Sigma, t)_i\delta_{v_i}, 
    \]   
    where
    \begin{align*} 
    \MA_\theta(\Sigma,t)_i&:=\frac{n!}{\norm[\alpha_i]\vol(\theta)}\int_\bbr \vol_{\bbr^{n-1}}(\Delta^\lambda(\nu_i;\theta;\Sigma,t)\cap \{x\in \R^n\mid \langle\alpha_i, x\rangle +t_i=\lambda\})\,d\lambda\\
    &=\frac{n!}{\vol(\theta)}\vol_{\bbr^n}\left(\left\{p\in\bbr^n\mid (\theta,\langle \alpha_i,p\rangle+t_i,t,p)\in \Delta(\nu_i;\Sigma)\right\}\right).
    \end{align*}
\end{definition}
The equality in the above definition follows again from the coarea formula applied to $p\mapsto (\theta,\langle \alpha_i,\pi(p)\rangle+t_i,t,p)$. 
As before, $\MA_\theta(\Sigma, t)$ defines a vector in $\R^d$, which we still denote by $\MA_\theta(\Sigma,t)$.
We now show that the Monge--Amp\`ere operator is well-defined, that is, $\MA_\theta(\Sigma, t)$ is a probability measure. This follows from the following extension of Theorem~\ref{thm:diff} to all of the big cone.
\begin{theorem}\label{thm:realma}
    For a fixed $\Sigma=\{v_i\}_{i=1}^d\subset\Xqm$ and $\theta\in \Bigc(X)$, we have $\vol(\theta; \Sigma, \cdot)\in C^1(\R^d)$ with $\nabla\vol(\theta; \Sigma, t)=\MA_\theta(\Sigma, t)\in \sigma_\Sigma$, where $\sigma_\Sigma\subset\R^d$ is the standard simplex. Further, if $\theta_k\to \theta$ in $\Bigc(X)$, then $\MA_{\theta_k}(\Sigma, t)\to \MA_\theta(\Sigma, t)$ weakly.
\end{theorem}
\begin{proof}
    Write $K_i(\theta,t):=\left\{p\in\bbr^n\mid (\theta,\langle \alpha_i,p\rangle+t_i,t,p)\in \Delta(\nu_i;\Sigma)\right\}$ throughout. We first note that the expression $\frac{n!}{\vol(\theta)}\vol_{\bbr^n}\left(K_i(\theta,t)\right)$
    defining $\MA_\theta(\Sigma,t)_i$ is usc in $(\theta,t)\in Big(X)_\bbr\times \bbr^{|\Sigma|}$ by standard convex body theory. Let $(\theta_k,t_k)\to (\theta,t)$. Since all the $K_i(\theta_k,t_k)$ and $K_i(\theta,t)$ are contained in a same compact set of $\bbr^n$ by Proposition \ref{prop:bigcvx} (as can be seen by normalising by $|\theta_k|$), Blaschke's selection theorem ensures existence of a subsequential limit $K_i'$ of the $K_i(\theta_k,t_k)$ under Hausdorff convergence. We then have $K_i'\subseteq K_i(\theta,t)$ by relative closedness of $\Delta(\nu_i;\Sigma)$.
    By continuity of the volume under Hausdorff convergence, we get
    $$\limsup_k \vol(K_i(\theta_k,t_k))=\vol(K_i')\leq \vol(K(\theta,t))$$
    as desired. 

    If $\theta-\theta'$ is ample, it is clear 
    that $K(\theta,t)\supseteq K(\theta',t)$, hence that
    $$\vol(\theta)\MA_\theta(\Sigma,t)_i\geq \vol(\theta')\MA_{\theta'}(\Sigma,t)_i.$$
    
    Let $\theta$ then be a big real class, and let $\theta_k\to \theta$ where $\theta_k\in \Bigc(X)_\bbq$ and $\theta_k-\theta$ is ample. Taking the limit in the above inequality gives
    $$\liminf_k \MA_{\theta_k}(\Sigma,t)_i \geq \MA_{\theta}(\Sigma,t)_i$$
    with upper semicontinuity giving the reverse inequality; in particular,
    $$\lim_k \MA_{\theta_k}(\Sigma,t)_i = \MA_\theta(\Sigma,t)$$
    for this decreasing approximation. By taking this limit in the probability measure equality one also obtains that
    $$\sum_i \MA_\theta(\Sigma,t)_i=1.$$
    In particular for each $i$, $\MA_\theta(\Sigma,t)_i=1-\sum_{j\neq i}\MA_\theta(\Sigma,t)_j$, hence it is also lsc, thus continuous. We conclude that $\vol(\theta; \Sigma, t)\in C^1$ as in \cite[Remark 4.29]{lazmus}.
\end{proof}
We also have the following extension of Proposition~\ref{prop:propertiesMA}.
\begin{proposition}\label{prop:manorm}
    For any $\theta\in \Bigc(X)$, the Monge--Amp\`ere operator $\MA_\theta\colon \cD(\Xqm)\to \cM^\qm$ descends to the Monge--Amp\`ere operator 
    \[
    \MA_\theta\colon \cN^\qm(\theta)\to \cM^\qm
    \]
    on the Hausdorff quotient by $\dd_{1, \theta}$. Further $\mathrm{supp} (\MA_\theta(\Sigma, t))=\rees(\theta; \Sigma, t)$ for any $(\Sigma, t)\in\cD(\Xqm)$.
\end{proposition}
\begin{proof}
    We first show that if $\dd_{1, \theta}((\Sigma, t), (\Sigma', t'))=0$, then $\MA_\theta(\Sigma, t)=\MA_\theta(\Sigma', t')$. After replacing one of the data by $(\Sigma\cup \Sigma', t\wedge t')$, we may assume $\Sigma\subseteq \Sigma'$ and $t'\le t$ on $\Sigma$. We may then further replace $\Sigma$ by $\Sigma'$ and $t$ by $(t, C{\bf 1})$ for $C\gg0$, so that the volume stays unchanged, by~\eqref{eq:embedding} and Proposition~\ref{prop:volcts}. Thus we are reduced to the case when $\Sigma=\Sigma'$ and $t\ge t'$. The result now follows from the same proof of Proposition~\ref{prop:propertiesMA} and using Theorem~\ref{thm:realma}.  

    Now let $(\Sigma=\{v_i\}_{i=1}^d, t)\in\cD(\Xqm)$. To see that the support of $\MA(\Sigma, t)$ is the set of Rees valuations, it will suffice to show $v_i\notin\rees(\theta; \Sigma, t)$ iff $\MA(\Sigma, t)_i=0$. To this end, note that $v_i\notin \rees(\theta; \Sigma, t)$ iff 
    \begin{equation}\label{eq:volred}
    \vol(\theta; \Sigma, t)=\vol(\theta; \Sigma\setminus\{v_i\}, t\setminus\{t_i\}).
    \end{equation}
    Since $\vol(\theta; \Sigma,t)$ is increasing, concave and eventually constant in each variable, we have~\eqref{eq:volred} holds iff $\vol(\theta; \Sigma, t+ce_i)=\vol(\theta; \Sigma\setminus\{v_i\}, t\setminus\{t_i\})$ for all $c\ge 0$, where $e_i\in\R^d$ is the $i$-th coordinate vector. This is in turn equivalent to $\MA(\Sigma, t)_i=0$ as desired.
    
\end{proof}

\subsection{The Calabi--Yau theorem for quasimonomial measures}

We are now ready to prove the following version of the non-Archimedean Calabi--Yau Theorem.
\begin{theorem}\label{thm:cyqm}
For any $\theta\in\Bigc(X)$, the Monge--Amp\`ere operator
\[
\MA_\theta\colon \cD(\Xqm)\to \cM^\qm, (\Sigma, t)\mapsto \MA_\theta(\Sigma, t)
\]
is surjective, and factors through 
\[
\MA_\theta\colon \cN^\qm(\theta)\to \cM^\qm.
\]
\end{theorem}

To prove Theorem \ref{thm:cyqm}, we rely on the following variational principle:
\begin{proposition}\label{prop:maenergy}
    Let $\Sigma=\{v_i\}_{i=1}^d\subset\Xqm$ be a finite set and $\xi\in\sigma_\Sigma$. Then $t\in \R^d$ satisfies 
    $\MA_\theta(\Sigma,t)=\xi$
    if and only if it achieves the supremum in
    $$\vol^\vee(\theta; \Sigma, \xi)\=\sup_{t\in \R^d} \{\vol(\theta;\Sigma,t) - \langle \xi,t\rangle\}=\max_{t\in K(\theta; \Sigma)}\{\vol(\theta;\Sigma,t) - \langle \xi,t\rangle\},$$
    where $K(\theta; \Sigma)\=[0, \max_{v\in \Sigma} T(\theta; v)]^d\subseteq \R^d$.
\end{proposition}
\begin{proof}
    Suppose first $\MA_\theta(\Sigma,t)=\xi$ for some $t\in \R^d$. Then by Proposition~\ref{prop:concave} and Theorem~\ref{thm:realma}, we have
    $$\vol(\theta;\Sigma,s)\leq \vol(\theta;\Sigma,t)+\langle \nabla \vol(\theta;\Sigma,t),s-t\rangle = \vol(\theta; \Sigma, t)+\langle\xi, s-t\rangle,$$
    for any $s\in \bbr^{d}$. This shows $\vol^\vee(\theta; \Sigma, \xi)=\vol(\theta; \Sigma, t)-\langle \xi, t\rangle$ after rearranging.

    Conversely, if $t$ computes $\vol^\vee(\theta; \Sigma, \xi)$, then by Theorem~\ref{thm:realma}, it maximizes the $C^1$ function $f(t)\=\vol(\theta; \Sigma, t)-\langle \xi, t\rangle$; in particular $f'(t)=0$, and this gives $\MA_\theta(\Sigma, t)=\xi$.

    Finally, it follows from translation invariance and~\eqref{eq:volconst} that the supremum is always attained on the compact set $K(\theta; \Sigma)$.
\end{proof}

\begin{proof}[Proof of Theorem~\ref{thm:cyqm}.] 
This follows from Proposition~\ref{prop:maenergy} and Proposition~\ref{prop:manorm}.
\end{proof}

We next turn to the uniqueness part of the Calabi--Yau theorem. We start by noting the following comparison between $\dd_{1, \theta}$-pseudometrics and the usual $\ell^1$ norm on the Euclidean space.
\begin{proposition}\label{prop:normcomparison}
    Fix a finite set $\Sigma\subset X^\qm$ of size $d$ and $\theta\in \Bigc(X)$. Then for any $t,t'\in \bbr^{d}$,
    \[
    \dd_{1,\theta}( (\Sigma,t),(\Sigma,t'))\le \|t'-t\|_1,
    \]
    where $\|\cdot\|_1$ is the standard $\ell^1$-norm on $\R^{d}$.
    If in addition $\MA_\theta(\Sigma, t), \MA_\theta(\Sigma, t')\in \sigma_\varepsilon$ for some $\varepsilon>0$, where 
    \[
    \sigma_\varepsilon\=\{\xi\in \sigma_{\Sigma}\mid \min_i\xi_i\ge \varepsilon\},
    \]
    then 
    \[
    \varepsilon \|t'-t\|_1\le \dd_{1,\theta}( (\Sigma,t),(\Sigma,t')).
    \]
\end{proposition}
\begin{proof}
   By Proposition~\ref{prop:concave}, for any $s_1, s_2\in \R^d$, we have
    \[
    \langle \nabla\vol(\theta; \Sigma, s_1), s_1-s_2\rangle \le \vol(\theta; \Sigma, s_1)-\vol(\theta; \Sigma, s_2)\le \langle\nabla\vol(\theta; \Sigma, s_2), s_1-s_2\rangle.
    \]
    By setting $s_1=t,t'$ and $s_2=t\wedge t'$ we get the desired comparison:
    \begin{align*}
        \varepsilon\|t-t'\|_1
        &\le\langle \nabla\vol(\theta; \Sigma, t), t-t\wedge t'\rangle+\langle \nabla\vol(\theta; \Sigma, t'), t'-t\wedge t'\rangle\\
        &\le \dd_{1,\theta}( (\Sigma,t),(\Sigma,t'))
        \le \langle \nabla\vol(\theta; \Sigma, t\wedge t'), |t-t'|\rangle\le \|t-t'\|_1,        
    \end{align*}
    where $|t-t'|\=(|t_{i}-t'_{i}|)_i$.
\end{proof}


\begin{theorem}\label{thm:masol}
    For any $\theta\in \Bigc(X)$, the Monge--Amp\`ere operator defines a bijection
    \[
    \MA_\theta\colon \cN^\qm(\theta)/\R\to \cM^\qm.
    \]
   In particular, if $\mu\in \cM^\qm$ has support $\Sigma$ with size $d$, there exists a unique $t\in \R^d$ up to translation such that $\MA_\theta(\Sigma, t)=\mu$.
\end{theorem}
\begin{lemma}\label{lem:equality}
    Let $\theta\in \Bigc(X)$. Suppose $\MA_\theta(\Sigma,t)=\MA_\theta(\Sigma,t')$. Then
    \[
    \alpha\mapsto \vol(\theta; \Sigma, (1-\alpha)t+\alpha t')
    \]
    is affine on $[0,1]$.
\end{lemma}
\begin{proof}
    This again follows immediately from Proposition~\ref{prop:concave}. Indeed, 
    \[
    \langle\nabla\vol(\theta; \Sigma, t'), t'-t\rangle\le \vol(\theta; \Sigma, t')-\vol(\theta; \Sigma, t)
    \le \langle\nabla\vol(\theta; \Sigma, t), t'-t\rangle.
    \]
    Since $\nabla\vol(\theta; \Sigma, t)=\MA_\theta(\Sigma, t)=\MA_\theta(\Sigma,t')=\nabla\vol(\theta; \Sigma, t')$, we must have
    \[
    \vol(\theta; \Sigma, t')-\vol(\theta; \Sigma, t)
    = \langle\nabla\vol(\theta; \Sigma, t), t'-t\rangle.
    \]
\end{proof}
\begin{proof}[Proof of Theorem~\ref{thm:masol}]
    Suppose $(\Sigma', t')\in \cD(\Xqm)$ is any other solution to $\MA_\theta(\Sigma', t')=\mu$ with $\vol(\theta; \Sigma', t')=\vol(\theta; \Sigma, t)$. Then by Proposition~\ref{prop:manorm} we have $\rees(\theta; \Sigma', t')=\Sigma$ and so $\Sigma'$ contains $\Sigma$ as a subset. Thus there exists $(\Sigma', t'')\in\cD(\Xqm)$ with $t''_v=t_v$ for all $v\in \Sigma$ with $\dd_{1, \theta}((\Sigma, t), (\Sigma', t''))=0$. Then by the triangle inequality we have
    \[
    \dd_{1,\theta}((\Sigma, t), (\Sigma', t'))\le \dd_{1,\theta}((\Sigma, t), (\Sigma', t''))+\dd_{1,\theta}((\Sigma', t''), (\Sigma', t'))=0,
    \]
    where the last equality follows from Lemma~\ref{lem:equality} and Proposition~\ref{prop:d1uniqueness}. This shows $(\Sigma, t)$ and $(\Sigma', t')$ define the same element in $\cN^\qm(\theta)/\R$.

    To prove the second assertion, assume now $\mu=\sum_{i=1}^d\xi_i\delta_{v_i}$ with $\xi_i>0$ for all $i$, so that $\mathrm{supp}(\mu)=\Sigma\=\{v_i\}_{i=1}^d$. Suppose $t_1, t_2\in \R^d$ are two solutions such that $\MA_\theta(\Sigma,t_1)=\MA_\theta(\Sigma, t_2)=\mu$.     
    After translation, we may assume $\vol(\theta; \Sigma, t_1)=\vol(\theta; \Sigma, t_2)$. Put also $\varepsilon=\min_i \xi_i>0$. Then by Proposition~\ref{prop:normcomparison} we get
    $$\varepsilon \norm[t_1-t_2]_1\leq \dd_{1,\theta}((\Sigma,t_1),(\Sigma,t_2))=0,$$
    and so $t_1=t_2$. 
\end{proof}

\begin{corollary}
    Let $\Sigma\subset\Xqm$ be finite. Then 
    \[
    \MA_\theta\colon (\cN^\Sigma(\theta)/\R, \qd_1)\to \sigma_\Sigma
    \]
    is a homeomorphism, where $\sigma_\Sigma$ is equipped with the Euclidean topology.
\end{corollary}
\begin{proof}
    By Theorem~\ref{thm:masol} and Proposition~\ref{prop:normcomparison}, we have that $\MA^{-1}_\theta\colon\sigma_\Sigma\to \cN^\Sigma(\theta)/\R$ is a continuous bijection. Thus this is a homeomorphism since $\sigma_\Sigma$ is compact and $\cN^\Sigma(\theta)/\R$ is Hausdorff.
\end{proof}


\begin{corollary}\label{coro:continuitysolutions}
    Let $(\mu_k)_k$ be a sequence of quasimonomial measures with support $\Sigma=\{v_i\}_{i=1}^d$ converging weakly to a quasimonomial measure $\mu$ with the same support. Assume $t_k\in \R^d$ solves $\MA_{\theta_k}(\Sigma, t_k)=\mu_k$, with $\theta_k\in \Bigc(X)$ converging to a class $\theta\in \Bigc(X)$. Then there exists a unique $t\in \R^d$ up to translation such that
    $$\MA_\theta(\Sigma,t)=\mu,\quad \qd_{1,\theta}((\Sigma,t_k),(\Sigma,t))\to 0.$$
\end{corollary}
\begin{proof}
    After translation, we may and will assume $\min_{i}t_{k,i}=0$. Note first that for $k\gg0$, we have $t_k\in[0, C]^d$ for some $C\gg0$. Indeed, for each $k$, since $\mathrm{supp}(\mu_k)=\Sigma$, we have $t_k\in[0, \max_iT(\theta_k, v_i)]^d$ by~\eqref{eq:volconst}, and we may pick $C=2\max_iT(\theta, v_i)$ by Lemma~\ref{prop:Tcts}. Thus by compactness, after passing to a subsequence, we may assume $t_k\to t$. It follows immediately from Proposition~\ref{prop:normcomparison} that
    $$\dd_{1,\theta}((\Sigma,t_k),(\Sigma,t))\leq \norm[t_k-t]_1\to 0.$$ 
    In view of Theorem~\ref{thm:masol}, it remains to show $\MA_\theta(\Sigma, t)=\mu$.

    Denote by $\xi_k\in\sigma_\Sigma$ the corresponding vectors defined by $\mu_k$.
    Since $\MA_{\theta_k}(\Sigma, t_k)=\mu_k$, by Proposition~\ref{prop:maenergy} we have for any $s\in \R^d$,
    \begin{equation}\label{eq:volmax}
        \vol(\theta_k;\Sigma,s)-\langle \xi_k,s\rangle \leq \vol(\theta_k;\Sigma,t_k)-\langle \xi_k,t_k\rangle.
    \end{equation}
    On the other hand, by concavity of the volume function, we have
    \begin{align*}
    \vol(\theta_k; \Sigma, t_k)-\vol(\theta; \Sigma, t)&=\vol(\theta_k; \Sigma, t_k)-\vol(\theta_k; \Sigma, t)+\vol(\theta_k; \Sigma, t)-\vol(\theta; \Sigma, t)\\
    &\le\langle\MA_{\theta_k}(\Sigma, t), t_k-t\rangle +\vol(\theta_k; \Sigma, t)-\vol(\theta; \Sigma, t).
    \end{align*}
    Thus $\vol(\theta_k; \Sigma, t_k)\to \vol(\theta; \Sigma, t)$ as $k\to \infty$ by Proposition~\ref{prop:manorm} and Proposition~\ref{prop:volcts}, and taking $k\to \infty$ in~\eqref{eq:volmax} gives    
    $\MA_\theta(\Sigma,t)=\mu$ by Proposition~\ref{prop:maenergy} again. 
\end{proof}

\subsection{Comparison with other definitions of Monge--Ampère measures}
Finally we end with showing that our notion of the Monge--Amp\`ere operator generalizes the existing ones in the literature. We start with the case when the polarization is ample.

\begin{proposition}
    If $\theta=c_1(L)$ with $L$ ample, then the Monge--Amp\`ere operator 
    \[
    \MA_\theta\colon \cN^\qm(L)\to \cM^\qm
    \]
    agrees with the one in~\cite{bj:kstab1}. 
\end{proposition}
\begin{proof}
    Let $\mu\=\sum_{i=1}^d\xi_i\delta_{v_i}\in \cM^\Sigma$. By Proposition~\ref{prop:maenergy} and~\cite[Proposition~7.8]{bj:kstab1}, it suffices to show that
    \[
    E^\vee(\mu)=\sup_{t\in \R^d}\{\vol(\theta; \Sigma, t)-\langle\xi, t\rangle\}.
    \]
    For any $\chi\in \cN_\R$, let $\chi'\=\min_{v_i\in \Sigma}\{v_i+\FS(\chi)(v_i)\}\in \cN^\Sigma(L)$. By~\cite[Theorem~2.16]{bj:kstab1}, we have $\chi'\ge \chi^\hom\ge \chi$ and so $\vol(\chi')\ge \vol(\chi)$. On the other hand, Proposition~6.5 in \emph{loc.cit.} shows that $\chi'=\min_{v_i\in \Sigma}\{v_i+\FS(\chi')(v_i)\}$ with $\FS(\chi')\le \FS(\chi)$ on $\Sigma$. Taking $t\=(\FS(\chi')(v_i))_{i=1}^d$ yields
    \[
    \vol(\theta; \Sigma, t)-\langle\xi, t\rangle = \vol(\chi')-\int \FS(\chi')\,\mu\ge  \vol(\chi)-\int \FS(\chi)\,\mu.
    \]
    Thus $E^\vee(\mu)\le \sup_{t\in\R^d}\{\vol(\theta; \Sigma, t)-\langle\xi, t\rangle\}$. Conversely, for any $t=(t_i)_{i=1}^d\in \R^d$, let $\chi=\chi_{\Sigma, t}$. Then by Proposition~\ref{prop:reesnorm} we have $\FS(\chi)(v_i)\le t_i$, and so
    \[
    \vol(\theta; \Sigma, t)-\langle \xi, t\rangle\le \vol(\chi)-\int\FS(\chi)\, \mu.
    \]
    This proves the reverse inequality.
\end{proof}

Let us also remark that, in the divisorial case, our notion agrees with the Monge--Ampère measure of test configurations for big line bundles as developed in \cite{derreb:3}. We refer the reader to the aforementioned references for further details on test configurations and their Monge--Ampère measures. Suppose $L$ is big and $X$ is normal, and $(\cX,\cL)$ is a test configuration in the sense of \cite[Definition 3.4]{derreb:3}, the \textit{non-Archimedean Monge--Ampère energy} is
$$E^\na_L(\cX,\cL):=\frac{\vol(\cL)}{(n+1)\vol(L)}.$$
Furthermore, denote by $\cX_0=\sum_{i=1}^d b_i E_i$ the central fiber of $\cX$ and $v_i\in X^\div$ the associated Rees valuations in the sense of \cite[Corollary 4.10]{bhj:duistermaat}, the Monge--Ampère measure is given by
$$\MA_L(\cX,\cL):=\vol(L)\mi\sum_i b_i \vol_{\cX|E_i}(\cL)\delta_{v_i}.$$
A test configuration also induces a divisorial norm $\chi_{\cX,\cL}=\chi_{\Sigma,t}$, given by the divisorial datum 
$$\left(\Sigma=(v_i)_i,\,t=(b_i\mi \ord_{E_i}(\cL-\pi^*L))_i\right).$$
Since $\MA_L(\cX,\cL)$ from \cite{derreb:3} is also the gradient of $E^\na_L(\cX,\pi^*L+\sum_i t_i b_i E_i)$ which agrees with $\vol(L;\Sigma,t)$ locally in $\bbr^{d}$ by~\cite[Lemma 4.12]{derreb:3}, one then has that
\begin{equation}
    \MA_L(\cX,\cL)=\MA_L(\Sigma,t).
\end{equation}

\bibliographystyle{alpha}
\bibliography{bib}

\end{document}